\pdfoutput=1
\RequirePackage{silence}
\documentclass[a4paper,svgnames]{amsart}
\usepackage[hmarginratio={1:1},vmarginratio={1:1},lmargin=60.0pt,tmargin=59.0pt]{geometry}

\usepackage{booktabs}
\allowdisplaybreaks

\usepackage[T1]{fontenc}
\usepackage{latexsym,exscale,amsfonts,amssymb,mathtools}
\usepackage{amsmath,amsthm,amsfonts,amssymb,amscd,textcomp,bbm}
\usepackage{mathrsfs,stackrel}
\usepackage{etoolbox}
\usepackage{tikz,tikz-cd}
\usetikzlibrary{matrix,arrows.meta,decorations.markings,shapes.multipart,decorations.pathreplacing,decorations.pathmorphing,shapes.geometric}
\usepackage{aliascnt}
\usepackage{ytableau}
\usepackage{xparse}
\usepackage{cite}

\usepackage{dynkin-diagrams}
\usepackage{nicematrix}
\usepackage{tikz-3dplot}
\usepackage{mathrsfs}
\DeclareMathAlphabet{\mathscrbf}{OMS}{mdugm}{b}{n}
\usepackage{float}
\usepackage{caption}

\usepackage{array}
\newcolumntype{C}{>{$}c<{$}}
\newcolumntype{P}[1]{>{\centering\arraybackslash}p{#1}}

\usepackage{xcolor,colortbl}

\definecolor{lightgreen}{HTML}{CCFFCC}
\definecolor{lightblue}{HTML}{CCCCFF}

\definecolor{mygray}{gray}{0.6}
\definecolor{mygraydark}{gray}{0.4}
\definecolor{mygraylight}{gray}{0.85}
\definecolor{spinach}{RGB}{46,139,87}
\definecolor{tomato}{RGB}{255,99,71}
\definecolor{orchid}{RGB}{143,40,194}
\definecolor{neon}{RGB}{77,77,255}
\definecolor{lightneon}{RGB}{110,110,255}
\definecolor{pumpkin}{RGB}{224,180,80}
\definecolor{citron}{RGB}{190,180,90}

\definecolor{lava}{RGB}{207,16,32}
\definecolor{cream}{RGB}{255,253,208}
\definecolor{verdigris}{RGB}{67,179,174}
\definecolor{Black}{RGB}{0,0,0}
\definecolor{mydarkblue}{RGB}{10,10,170}
\definecolor{darkspinach}{RGB}{20,70,20}
\definecolor{darktomato}{RGB}{155,40,30}
\definecolor{darkorchid}{RGB}{50,10,100}
\definecolor{darklava}{RGB}{150,8,16}

\definecolor{zero}{RGB}{0,0,0}
\definecolor{one}{RGB}{255,0,0}
\definecolor{two}{RGB}{0,255,0}
\definecolor{three}{RGB}{0,0,255}

\definecolor{lightergray}{gray}{0.9}

\usepackage{todonotes}

\usepackage{enumitem}
\setlist[enumerate]{itemsep=0.15cm,label=\emph{\upshape(\alph*)}}
\setlist[enumerate,2]{itemsep=0.15cm,label=\emph{\upshape(\roman*)}}
\setlist[enumerate,3]{itemsep=0.15cm,label=\emph{\upshape(\Alph*)}}

\let\emph\relax
\DeclareTextFontCommand{\emph}{\bfseries\em}

\renewcommand{\dots}{\text{...}}

\newcommand{\placeholder}{{}_{-}}
\newcommand{\mystrut}{\rule[-0.2\baselineskip]{0pt}{0.9\baselineskip}}

\usepackage[all]{xy}
\usepackage{tikz}
\usetikzlibrary{cd}
\usetikzlibrary{decorations}
\usetikzlibrary{decorations.markings}
\usetikzlibrary{decorations.pathreplacing}
\usetikzlibrary{decorations.pathmorphing}
\usetikzlibrary{arrows.meta,shapes,positioning,matrix,calc}
\usetikzlibrary{shapes.callouts}
\tikzset{anchorbase/.style={baseline={([yshift=-0.5ex]current bounding box.center)}},
tinynodes/.style={font=\tiny,text height=0.25ex,text depth=0.05ex},
smallnodes/.style={font=\scriptsize,text height=0.75ex,text depth=0.15ex},
usual/.style={line width=2.0,color=black},
crossline/.style={preaction={draw=white,line width=5.75pt,-}},
}
\tikzstyle directed=[postaction={decorate,decoration={markings,mark=at position #1 with {\arrow[line width=0.3mm, black]{>}}}}]

\newcommand{\R}{\mathbb{R}}

\newcommand{\N}{\mathbb{Z}_{\geq 0}}

\newcommand{\Z}{\mathbb{Z}}

\newcommand{\Ff}{\mathbb{F}}

\def\NewTheorem#1{%
\newaliascnt{#1}{equation}%
\newtheorem{#1}[#1]{#1}%
\aliascntresetthe{#1}%
\expandafter\def\csname #1autorefname\endcsname{#1}%
}
\def\equationautorefname~#1\null{(#1)\null}

\numberwithin{equation}{subsection}

\NewTheorem{Proposition}
\NewTheorem{Theorem}
\NewTheorem{Corollary}
\AtEndEnvironment{Corollary}{\null\hfill$\square$}%
\NewTheorem{Lemma}
\NewTheorem{Conjecture}
\NewTheorem{Speculation}
\NewTheorem{Question}
\NewTheorem{Assumption}
\theoremstyle{definition}
\NewTheorem{Definition}
\AtEndEnvironment{Definition}{\null\hfill$\Diamond$}%
\NewTheorem{Classification Problem}
\AtEndEnvironment{Classification Problem}{\null\hfill$\Diamond$}%
\NewTheorem{Notation}
\AtEndEnvironment{Notation}{\null\hfill$\Diamond$}%
\NewTheorem{Example}
\AtEndEnvironment{Example}{\null\hfill$\Diamond$}%
\NewTheorem{Examples}
\AtEndEnvironment{Examples}{\vskip-10mm\null\hfill$\Diamond$}%

\theoremstyle{remark}
\NewTheorem{Remark}
\AtEndEnvironment{Remark}{\null\hfill$\Diamond$}%

\usepackage[hypertexnames=false]{hyperref}
\usepackage{bookmark}
\hypersetup{
pdftoolbar=true,
pdfmenubar=true,
pdffitwindow=false,
pdfstartview={FitH},
pdftitle={On detection probabilities of link invariants},
pdfauthor={Tuomas Kelom{\"a}ki, Abel Lacabanne, Daniel Tubbenhauer, Pedro Vaz and Victor L. Zhang},
pdfsubject={},
pdfcreator={Tuomas Kelom{\"a}ki, Abel Lacabanne, Daniel Tubbenhauer, Pedro Vaz and Victor L. Zhang},
pdfproducer={Tuomas Kelom{\"a}ki, Abel Lacabanne, Daniel Tubbenhauer, Pedro Vaz and Victor L. Zhang},
pdfkeywords={},
pdfnewwindow=true,
colorlinks=true,
linkcolor=mydarkblue,
citecolor=teal,
filecolor=magenta,
urlcolor=orchid,
linkbordercolor=lava,
citebordercolor=teal,
urlbordercolor=orchid,
linktocpage=true
}

\def\makeautorefname#1#2{\csdef{#1autorefname}{#2}}
\makeautorefname{section}{Section}%
\makeautorefname{subsection}{Section}%
\makeautorefname{subsubsection}{Section}%

\begin{document}
\title[On detection probabilities of link invariants]{On detection probabilities of link invariants}
\author[T. Kelom{\"a}ki, A. Lacabanne, D. Tubbenhauer, P. Vaz and V.L. Zhang]{Tuomas Kelom{\"a}ki, Abel Lacabanne, Daniel Tubbenhauer, Pedro Vaz and Victor L. Zhang}

% \address{A.L.: Laboratoire de Math{\'e}matiques Blaise Pascal (UMR 6620), Universit{\'e} Clermont Auvergne, Campus Universitaire des C{\'e}zeaux, 3 place Vasarely, 63178 Aubi{\`e}re Cedex, France, \href{http://www.normalesup.org/~lacabanne}{www.normalesup.org/$\sim$lacabanne},
% \href{https://orcid.org/0000-0001-8691-3270}{ORCID 0000-0001-8691-3270}}
% \email{abel.lacabanne@uca.fr}

% \address{T.K.: Aalto University, Department of Mathematics and Systems Analysis, Otakaari 1, 02150, Espoo, Finland}
% \email{tuomas.kelomaki@aalto.fi}

% \address{D.T.: The University of Sydney, School of Mathematics and Statistics F07, Office Carslaw 827, NSW 2006, Australia, \href{http://www.dtubbenhauer.com}{www.dtubbenhauer.com}, \href{https://orcid.org/0000-0001-7265-5047}{ORCID 0000-0001-7265-5047}}
% \email{daniel.tubbenhauer@sydney.edu.au}

% \address{P.V.: Institut de Recherche en Math{\'e}matique et Physique, 
% Universit{\'e} catholique de Louvain, Chemin du Cyclotron 2,  
% 1348 Louvain-la-Neuve, Belgium, \href{https://perso.uclouvain.be/pedro.vaz}{https://perso.uclouvain.be/pedro.vaz}, \href{https://orcid.org/0000-0001-9422-4707}{ORCID 0000-0001-9422-4707}}
% \email{pedro.vaz@uclouvain.be}

% \address{V.L.Z.: University of New South Wales (UNSW), School of Mathematics and Statistics, NSW 2052, Australia, \href{https://dustbringer.github.io/}{dustbringer.github.io} }
% \email{victor.l.zhang@unsw.edu.au}

\begin{abstract}
We prove that, for many standard link invariants, both the proportion of distinct invariant values and the detection probability among prime alternating links with at most $n$ crossings decay exponentially in $n$, with an explicit universal rate. In fact, almost every such link belongs to an invariant fiber whose size is itself exponential in $n$. This phenomenon applies broadly, in particular to the Jones and HOMFLYPT polynomials and integral Khovanov homology. The companion website gives a much more detailed view of the data, including complete distributions of fiber sizes, separate alternating and non-alternating data, and topological data analysis.
\end{abstract}


\addtocontents{toc}{\protect\setcounter{tocdepth}{1}}

\maketitle

\tableofcontents

\section{Introduction}\label{S:Intro}

Our goal in this paper is to understand the effectiveness of knot invariants (classical, quantum, and beyond) as tools for distinguishing knots. For a finite set $X$ and an invariant $Q$, write
\[
v_Q(X)=\frac{\#Q(X)}{\#X},\qquad
d_Q(X)=\frac{\#\big\{x\in X\mid Q^{-1}\big(Q(x)\big)\cap X=\{x\}\big\}}{\#X}.
\]
Thus $v_Q(X)$ is the \emph{value ratio}, while $d_Q(X)$ is the actual \emph{detection ratio}; in particular $d_Q(X)\leq v_Q(X)$.

We use the following terminology throughout. The invariant $Q$ \emph{distinguishes} $x,y\in X$ if $Q(x)\neq Q(y)$, and it \emph{detects} $x$ in $X$ if $Q^{-1}(Q(x))\cap X=\{x\}$. Thus detection means distinction from every other element of $X$. Viewed as a classifier, $Q$ partitions $X$ into the \emph{fibers} $Q^{-1}(q)\cap X$, for $q\in Q(X)$. Whenever we speak of a fiber, the invariant and ambient set are the ones under discussion.

The following theorem captures our main theoretical result; throughout, knots and links are unoriented, and links are non-split.

\begin{Theorem}\label{T:Main}
Set
\[
\rho=0.9943699249\ldots\quad (\rho^{-1}=1.00566195232\ldots),
\]
where the exact definition is given in \autoref{T:Decay}. For every invariant marked green in \autoref{Tab:Main} and every $\delta\in(\rho,1)$, both the value ratio and the detection ratio on prime alternating links with at most $n$ crossings are $O(\delta^n)$. Moreover, for every $\eta\in(1,\rho^{-1})$, the probability that a uniformly random such link belongs to a fiber of size smaller than $\eta^n$ is $O(\theta^n)$ for some $\theta\in(0,1)$.
\end{Theorem}

\begin{table}[!ht]
\centering
\begin{tabular}{cc}
\rowcolor{lightergray} Polynomial & P(detection$)=0$? \\
Jones & \cellcolor{lightgreen}Yes \\
Alexander & \cellcolor{lightgreen}Yes \\
$\mathrm{SL}_N$ & \cellcolor{lightgreen}Yes \\
HOMFLYPT & \cellcolor{lightgreen}Yes \\
Jones all colors & \cellcolor{lightgreen}Yes \\
$\mathrm{SL}_N$ for $(1^k)$ & \cellcolor{lightgreen}Yes \\
$\mathrm{SL}_N$ for $(k)$ & \cellcolor{lightgreen}Yes \\
HOMFLYPT for $(1^k)$ & \cellcolor{lightgreen}Yes \\
HOMFLYPT for $(k)$ & \cellcolor{lightgreen}Yes \\
Kauffman & \cellcolor{lightgreen}Yes \\
$\mathrm{SO}_{2N+1}$ & \cellcolor{lightgreen}Yes \\
$\mathrm{SP}_{2N}$ & \cellcolor{lightgreen}Yes \\
$\mathrm{SO}_{2N}$ & \cellcolor{lightgreen}Yes \\
$\mathrm{G}_2$ & \cellcolor{lightgreen}Yes \\
\end{tabular}
\begin{tabular}{cc}
\rowcolor{lightergray} Homology & P(detection$)=0$? \\
Khovanov over $\Z$ & \cellcolor{lightgreen}Yes \\
Odd Khovanov over $\Z$ & \cellcolor{lightgreen}Yes \\
$\widehat{\mathrm{HFK}}$ over $\Ff_2$ & \cellcolor{lightgreen}Yes \\
Khovanov--Rozansky over $\Z$ & \cellcolor{blue!70}Likely \\
HOMFLYPT over $\Z$ (for knots) & \cellcolor{blue!70}Likely \\
$(1^k)$ or $(k)$ versions of these & \cellcolor{blue!70}Likely \\
Khovanov--Floer theories over $\Ff_2$ & \cellcolor{blue!70}Likely \\
\hline
\rowcolor{lightergray} Other & P(detection$)=0$? \\
Signature & \cellcolor{lightgreen}Yes \\
Determinant & \cellcolor{lightgreen}Yes \\
Double branched cover & \cellcolor{lightgreen}Yes \\
HF of double branched cover & \cellcolor{lightgreen}Yes \\
Many invariants that we forgot to add & \cellcolor{cream!70}Depends, but likely \\
\end{tabular}
\\[0.1cm]
\caption{Summary of our main results. Green entries (``Yes'') are proved in this paper. Blue entries (``Likely'') indicate invariants for which we outline proof strategies and supply big data evidence, though final arguments remain open. Here and below, a random link is chosen uniformly from the prime alternating links with at most $n$ crossings, and then $n\to\infty$.}
\label{Tab:Main}
\end{table}

\subsection{What this means, roughly speaking}

Classical invariants such as the Alexander polynomial, signature and determinant have long been used to distinguish knots. The Jones polynomial \cite{Jo-jones-polynomial} introduced a different family coming from quantum groups and led to invariants such as the Witten--Reshetikhin--Turaev invariants \cite{ReTu-invariants-3-manifolds-qgroups}. Khovanov homology \cite{Kh-cat-jones} and related constructions categorify these polynomial invariants.

A natural question is:
\emph{How effective are such invariants (classical, quantum, and beyond) as tools for distinguishing knots?}
Distinguishing knots is certainly not the only reason to study these invariants, but it is a natural thing to measure.

A form of ``folk wisdom'' (with increasingly rigorous and empirical backing, see \cite{St-number-polynomials,DlGuSa-data,TuZh-knot-data} for a very biased sample) suggests:
\begin{gather*}
\fcolorbox{orchid!50}{spinach!10}{\mystrut``Many invariants detect a random knot with probability zero.''}
\end{gather*}
That is: put all knots in a bag, shake, and draw one. How likely is a given invariant (say, the determinant or Khovanov homology) to single it out? For the unknot and Khovanov homology the answer is famously positive \cite{KrMo-unknot-detector}, but the unknot is exceptional; a generic knot typically shares its invariants with many others. Moreover, they share their invariant with exponentially many other links.

We can prove this for prime alternating links. The full set of knots is intractable, but prime alternating links can be counted asymptotically \cite{Tu-planar,SuTh-growth}. We therefore restrict attention to alternating links, where our main theorem \autoref{T:Main} applies. The key tool is oriented mutation: following \cite{St-number-polynomials}, we show that the number of prime alternating links modulo oriented mutation has a strictly smaller exponential growth rate than the total number of prime alternating links. Consequently, invariants insensitive to oriented mutation have exponentially decaying value and detection ratios.

Alternating links and oriented mutation should be viewed as methodological devices rather than the ultimate focus. (We expect most of the invariants in this paper to also have vanishing value ratio on knots up to mutation, cf. \autoref{table:percent-jones-mut}.) To probe the situation for knots themselves, we complement the theory with large-scale ``big data'' computational experiments in the spirit of \cite{LeHaSa-jones,DlGuSa-mapper,LaTuVa-big-data,DlGuSa-data,TuZh-knot-data}. These experiments concentrate on quantum invariants, comparing polynomial invariants to their categorifications and measuring value and detection ratios across extensive datasets. The homological invariants add much less distinguishing power than their size and computational cost might suggest. For example, through 16 crossings rational Khovanov homology raises the Jones detection ratio by only $5.36$ percentage points, from $29.09\%$ to $34.45\%$. For further statistics, the companion website is worth opening: it contains the full distributions of fiber sizes, separate alternating and non-alternating data, mutation quotients, and several topological data analyses. The plots reproduced here show only a small part of it. This approach and its results are arguably more important than the ``folk theorem'' \autoref{T:Main} itself.

\subsection{Some remarks}

Before we prove \autoref{T:Main}, let us collect some remarks.

\begin{Remark}
Links are easier to study asymptotically than knots, but we do not expect a huge difference between their detection rates. The important restriction is to nonsplit alternating links. We do not know how to adapt the proof to non-alternating links; see \autoref{R:Knots}. We expect the detection rate to be worse there; see \autoref{S:BigData} and \cite{DlGuSa-data,TuZh-knot-data} for evidence from the data.
\end{Remark}

\begin{Remark}
Alternatively, one could ask whether a given invariant distinguishes two fixed knots and not between one knot and all other knots. See e.g. \cite{TuZh-knot-data} for some results in this direction.
\end{Remark}

\begin{Remark}\label{R:NotSurprising}
(\emph{Why this phenomenon is not surprising.})
There are several reasons to expect that many knot invariants perform poorly as classifiers. Recorded below is one concrete reason, while another is the ``folk wisdom'' mentioned above.

One would expect, say, quantum invariants to perform exceptionally well on knots with few crossings, and indeed they do (cf.\ \autoref{S:BigData}). In contrast, detection across the entire set of knots is widely expected to be very difficult, even if formal results are limited. As evidence, many nearby decision/optimization problems are already hard: for example, knot genus is NP-hard \cite{MR2219001}, and the (graph) crossing number admits no constant-factor approximation unless $\textup{P}=\textup{NP}$ \cite{MR3017916}. Such barriers make it unrealistic to expect generic invariants, many of which are ``efficiently'' computable in practice, to act as effective classifiers at scale. In contrast, strong invariants provide a practical remedy, e.g. the knot complement is a complete invariant of a knot by \cite{MR965210}, but the associated 3-manifold decision problems are expected to be very difficult, see, for example, \cite{MR3609909} for a summary.
\end{Remark}

\begin{Remark}\label{R:Surprising}
(\emph{Why this phenomenon is surprising.})
We had no a priori reason to expect the invariants in \autoref{T:Main} to have detection probability zero, so the theorem was genuinely surprising to us (even if there were hints). To emphasize why \autoref{R:NotSurprising} does not translate into statements about detection probabilities, we give an analogy: deciding whether a graph contains a Hamiltonian cycle is NP-complete, yet in the standard random-graph models the heuristic ``Every graph is Hamiltonian'' is asymptotically correct with probability one (i.e. asymptotically almost surely). 

The finite data already contain large collisions. The Jones value ratio is $41.62\%$ through 18 crossings. Through 16 crossings, where we have the full distribution of fiber sizes, the actual Jones detection ratio is only $29.09\%$. Moreover, the Jones fiber-size distribution is close to a straight line on log-log axes, so the collisions occur at every visible scale rather than in a small exceptional set; see \autoref{S:BigData}.
\end{Remark}

\begin{Remark}
The restriction on colors in \autoref{Tab:Main} matters. Symmetric and exterior colors give the multiplicity-free cases used below, whereas non-rectangular colors need not be mutation invariant. For example, the $(2,1)$-colored HOMFLYPT polynomial distinguishes the Kinoshita--Terasaka and Conway mutant pair; see \cite{MR1395780,NaRaSi-colored-mutants}.
\end{Remark}

%\begin{Remark}
%Our arguments, and therefore \autoref{T:Main}, also apply to other invariants such as hyperbolic volume, (most) homological link invariants which admit spectral sequences from Khovanov homology, or the Heegaard--Floer homology of the branched double-cover associated to the link. This follows from the proof of \autoref{T:Main} and \cite{MR906585,Saltz,MR3049933}.
%\end{Remark}

\noindent\textbf{Acknowledgments.}
We thank Radmila Sazdanovic, whose talk (York, July 2024) sparked this project. We used ChatGPT to assist with coding, proofreading, conjecture generation, and aspects of proof development, and we gratefully acknowledge this support. We thank Alexis Langlois-Rémillard for a pointer to the literature, which sparked part of this collaboration. Part of this project was done while the second author visited the fourth author at the Universit{\'e} catholique de Louvain, and their hospitality is gratefully acknowledged.

This work used the Katana computational cluster (DOI: \url{https://doi.org/10.26190/669X-A286}), supported by Research Technology Services at UNSW Sydney, and some of the computations were performed using the resources of the Laboratoire de Math{\'e}matiques Blaise Pascal.

DT acknowledges support from the ARC Future Fellowship FT230100489 and suspects their role will soon be outsourced to AI.

\section{Proof of the main result}\label{S:Main}

We now explain and prove the statement {announced in the introduction}, and even the 
generalization in \autoref{T:Main} of it. But before we start:

\begin{Remark}
This paper is largely self-contained; however, we assume the reader has some familiarity with knot theory and analytic combinatorics. For background in these areas, we recommend standard references such as \cite{Ad-knots,FlSe-analytic-combinatorics}. Beginning after \autoref{T:Decay}, we also rely on several standard tools from the theory of quantum invariants. For background on these, see for example \cite{Tu-qgroups-3mfds,Oh-quantum-invariants,Tu-qt} for the polynomial invariants and \cite{BaNa-categorification-jones,Kh-link-homology} for homological invariants.
\end{Remark}

Our proof is a variation of arguments from \cite{Tu-planar,SuTh-growth,St-number-polynomials}, where the main work was done. The new point in the counting argument is that the upper bound can still be obtained using oriented mutation. This is the version needed for the invariants in \autoref{Tab:Main}.

The argument has four ingredients. First, the flyping theorem of Menasco--Thistlethwaite \cite{MeTh-flyping} allows us to pass from reduced alternating diagrams to links by quotienting by flypes. Second, Sundberg--Thistlethwaite \cite{SuTh-growth} give a growth constant $\lambda>6.1479$ for prime alternating links. Third, a refinement of Stoimenow's count \cite{St-number-polynomials}, using all permutations at rational nodes, gives an upper growth constant $\mu<6.114$ modulo oriented mutation. Finally, an invariant that is constant on oriented mutation classes takes at most as many values as there are classes. The strict inequality $\mu<\lambda$ gives exponential decay. The details below recall the two counts and explain the oriented mutation step.

\begin{Lemma}\label{L:EasySandwich}
Fix two sequences of sets $X_n,Y_n$ with $n\in\N$, and let $\alpha,\beta\in\R_{\geq 1}$. Then we have
\begin{gather*}
\Big(\limsup_{n\to\infty}\sqrt[n]{\# Y_n}\leq
\alpha<\beta\leq\liminf_{n\to\infty}\sqrt[n]{\# X_n}\Big)
\text{ implies }
\Big(\frac{\# Y_n}{\# X_n}\in O(\gamma^n)\text{ for some }\gamma\in(0,1)\Big).
\end{gather*}
Moreover, $\gamma$ can be chosen to be $\gamma=\frac{\alpha + \varepsilon}{\beta - \varepsilon}$ for some $\varepsilon>0$.

The same holds when replacing $\#Y_n$ and $\#X_n$ by $\sum_{k=0}^n\#Y_k$ and 
$\sum_{k=0}^n\#X_k$ on the right-hand side.
\end{Lemma}

\begin{proof}
%
%Let \( a_n = \#Y_n \), \( b_n = \#X_n \), and choose \( \varepsilon > 0 \) such that \( \alpha + \varepsilon < \beta - \varepsilon \). Set \( \gamma = \frac{\alpha + \varepsilon}{\beta - \varepsilon} \in (0,1) \). By the assumptions, \( \limsup_{n \to \infty} a_n^{1/n} \leq \alpha < \beta \leq \liminf_{n \to \infty} b_n^{1/n} \), so for sufficiently large \( n \), we have \( a_n^{1/n} \leq \alpha + \varepsilon \) and \( b_n^{1/n} \geq \beta - \varepsilon \). Hence, \( \frac{a_n}{b_n} \leq \gamma^n \), so \( \#Y_n / \#X_n \in O(\gamma^n) \). For the final claim one applies the geometric series.

Choose $\varepsilon>0$ with $\alpha+\varepsilon<\beta-\varepsilon$ and set
$r_1:=\alpha+\varepsilon$, $r_2:=\beta-\varepsilon$; then $1<r_1<r_2$. 
By the $\limsup/\liminf$ hypotheses, for $n\gg 0$ we have 
$\#Y_n\le r_1^{\,n}$ and $\#X_n\ge r_2^{\,n}$, hence
\[
\frac{\#Y_n}{\#X_n}\le \Bigl(\frac{r_1}{r_2}\Bigr)^{\!n}=:\gamma^n,
\qquad \gamma\in(0,1).
\]
For the summed version, note that
\[
\sum_{k=0}^n \#Y_k \le \sum_{k=0}^n r_1^{\,k}
= \frac{r_1^{\,n+1}-1}{r_1-1}\le \frac{r_1}{r_1-1}\, r_1^{\,n},
\qquad
\sum_{k=0}^n \#X_k \ge \#X_n \ge r_2^{\,n}.
\]
Therefore
\[
\frac{\sum_{k=0}^n \#Y_k}{\sum_{k=0}^n \#X_k}
\le \frac{r_1}{r_1-1}\Bigl(\frac{r_1}{r_2}\Bigr)^{\!n}
= C\,\gamma^n
\]
with the same $\gamma\in(0,1)$ and a constant $C>0$ independent of $n$.
\end{proof}

Let \((X_n)_{n\in\N}\) be a sequence of nonempty finite sets such that $X_{n}\subset X_{n+1}$ for all $n\in\N$.
Define the limiting set
\[
X := \bigcup_{n=0}^\infty X_n.
\]
Let $Y$ be a subset of $X$ and define, for all $n\in \N$, $Y_{n} = Y\cap X_{n}$.

\begin{Lemma}\label{L:EasySandwich2}
Suppose the hypotheses of \autoref{L:EasySandwich} hold for the nested finite sets $Y_n\subseteq X_n$. Then the probability that a uniformly random $x\in X_n$ lies in $Y_n$ tends to zero exponentially.
\end{Lemma}

\begin{proof}
We have that the relative sizes decay exponentially, i.e. for large enough $n$:
\[
\frac{\# Y_n}{\# X_n} \leq C\cdot \gamma^n, \quad \text{for some constants } C\in\R_{>0},\gamma\in(0,1).
\]
Then \(Y\) has asymptotic density zero in \(X\), in the sense that
\[
\lim_{n \to \infty} \frac{\#(Y \cap X_n)}{\# X_n} = \lim_{n \to \infty} \frac{\# Y_n}{\# X_n} = 0.
\]
In particular, if we equip \(X_n\) with the uniform probability measure, meaning
\[
\mathbb{P}_n(x) = 1/\#X_n,\quad x\in X_n,
\]
then the probability that a randomly chosen element of \(X_n\) lies in \(Y_n\) tends to zero exponentially:
\[
\mathbb{P}_n(x\in Y) = \mathbb{P}_n(x\in Y_n) = \frac{\#Y_n}{\#X_n} \leq C\cdot\gamma^n.
\]
In probabilistic terms, this means:
\[
\lim_{n\to\infty}\mathbb{P}_n(x\in Y) = 0.
\]
{This is the only notion of probability used in the paper: for each $n$, we choose uniformly from the finite set $X_n$, and then let $n\to\infty$.} This proves the claim.
\end{proof}

Below, we will only consider tangles with four boundary points and call them \emph{tangles} for short. We draw them in circles where one boundary string is marked (so that we can number the boundary strings, say clockwise starting from the marked one). Here is an example, the tangle is called $T$:
\begin{gather*}
\begin{tikzpicture}[anchorbase,scale=1]
\draw[usual,tomato] (0.25,0.4) to (0.25,0.6);
\draw[usual] (-0.25,0.4) to (-0.25,0.6);
\draw[usual] (0.25,-0.6) to (0.25,0.1);
\draw[usual] (-0.25,-0.6) to (-0.25,0.1);
\draw[fill=cream] (0,0) circle (0.5cm);
\node at (0,0) {$T$};
\end{tikzpicture}
=
\begin{tikzpicture}[anchorbase,scale=1]
\draw[usual,tomato] (0.4,0.25) to (0.6,0.25);
\draw[usual] (0.4,-0.25) to (0.6,-0.25);
\draw[usual] (-0.6,0.25) to (0.1,0.25);
\draw[usual] (-0.6,-0.25) to (0.1,-0.25);
\draw[fill=cream] (0,0) circle (0.5cm);
\node at (0,0) {$T$};
\end{tikzpicture}
=
\begin{tikzpicture}[anchorbase,scale=1]
\draw[usual] (0.5,-0.5) to (-0.5,0.5);
\draw[usual] (-0.5,-0.5) to (0,0);
\draw[usual,tomato] (0.5,0.5) to (0,0);
\draw[fill=cream] (0,0) circle (0.5cm);
\node at (0,0) {$T$};
\end{tikzpicture}
.
\end{gather*}
The building blocks that we use are
\begin{gather*}
0\colon
\begin{tikzpicture}[anchorbase,scale=1]
\draw[usual] (0.5,0) to[out=90,in=270] (0.5,0.5);
\draw[usual,crossline] (0,0) to[out=90,in=270] (0,0.5);
\end{tikzpicture}
,\quad
\Big(\infty\colon
\begin{tikzpicture}[anchorbase,scale=1]
\draw[usual] (0.5,0.5) to (0,0.5);
\draw[usual,crossline] (0,0) to (0.5,0);
\end{tikzpicture}
,\Big)\quad
1\colon
\begin{tikzpicture}[anchorbase,scale=1]
\draw[usual] (0.5,0) to[out=90,in=270] (0,0.5);
\draw[usual,crossline] (0,0) to[out=90,in=270] (0.5,0.5);
\end{tikzpicture}
,\quad
-1\colon
\begin{tikzpicture}[anchorbase,scale=1]
\draw[usual] (0,0) to[out=90,in=270] (0.5,0.5);
\draw[usual,crossline] (0.5,0) to[out=90,in=270] (0,0.5);
\end{tikzpicture}
.
\end{gather*}
The basic operations on such tangles are horizontal and vertical concatenation, called sum, defined as follows. For tangles $S$ and $T$ (where we use dashed lines for the marked strand of $S$): 
\begin{gather}\label{Eq:Sums}
S+_hT=
\begin{tikzpicture}[anchorbase,scale=1]
\draw[usual] (-0.1,0) to (0.6,0);
\draw[usual,densely dashed] (-0.1,0.5) to (0.6,0.5);
\draw[usual] (0.9,0) to (1.6,0);
\draw[usual,tomato] (0.9,0.5) to (1.6,0.5);
\draw[usual] (-1.1,0) to (-0.4,0);
\draw[usual] (-1.1,0.5) to (-0.4,0.5);
\draw[fill=cream] (-0.5,0.25) circle (0.5cm);
\draw[fill=cream] (1.0,0.25) circle (0.5cm);
\node at (-0.5,0.25) {$S$};
\node at (1.0,0.25) {$T$};
\end{tikzpicture}
,\quad
S+_vT=
\begin{tikzpicture}[anchorbase,scale=1]
\draw[usual,densely dashed] (0.25,0.4) to (0.25,1.1);
\draw[usual] (-0.25,0.4) to (-0.25,1.1);
\draw[usual,tomato] (0.25,1.4) to (0.25,2.1);
\draw[usual] (-0.25,1.4) to (-0.25,2.1);
\draw[usual] (0.25,-0.6) to (0.25,0.1);
\draw[usual] (-0.25,-0.6) to (-0.25,0.1);
\draw[fill=cream] (0,0) circle (0.5cm);
\draw[fill=cream] (0,1.5) circle (0.5cm);
\node at (0,0) {$S$};
\node at (0,1.5) {$T$};
\end{tikzpicture}
.
\end{gather}
Recall that a \emph{rational tangle} is a tangle obtained by (horizontally or vertically) adding the tangles $\pm 1$ to the tangle $0$. Such a tangle is \emph{prime} if it cannot be composed into two different nontrivial (neither $0$ nor $\infty$) summands using \autoref{Eq:Sums}.

Let us use the following notation for rational tangles 
(illustrated with three crossings; the evident generalization applies):
\begin{gather*}
3\colon
\raisebox{-0.25cm}{\rotatebox{90}{\begin{tikzpicture}[anchorbase,scale=1]
\draw[usual] (0.5,0) to[out=90,in=270] (0,0.5);
\draw[usual,crossline] (0,0) to[out=90,in=270] (0.5,0.5);
\draw[usual] (0.5,0.5) to[out=90,in=270] (0,1);
\draw[usual,crossline] (0,0.5) to[out=90,in=270] (0.5,1);
\draw[usual] (0.5,1) to[out=90,in=270] (0,1.5);
\draw[usual,crossline] (0,1) to[out=90,in=270] (0.5,1.5);
\end{tikzpicture}}}
,\quad
-3\colon
\raisebox{-0.25cm}{\rotatebox{90}{\begin{tikzpicture}[anchorbase,scale=1]
\draw[usual] (0,0) to[out=90,in=270] (0.5,0.5);
\draw[usual,crossline] (0.5,0) to[out=90,in=270] (0,0.5);
\draw[usual] (0,0.5) to[out=90,in=270] (0.5,1);
\draw[usual,crossline] (0.5,0.5) to[out=90,in=270] (0,1);
\draw[usual] (0,1) to[out=90,in=270] (0.5,1.5);
\draw[usual,crossline] (0.5,1) to[out=90,in=270] (0,1.5);
\end{tikzpicture}}}
,\quad
\bar{3}\colon
\begin{tikzpicture}[anchorbase,scale=1]
\draw[usual] (0,0) to[out=90,in=270] (0.5,0.5);
\draw[usual,crossline] (0.5,0) to[out=90,in=270] (0,0.5);
\draw[usual] (0,0.5) to[out=90,in=270] (0.5,1);
\draw[usual,crossline] (0.5,0.5) to[out=90,in=270] (0,1);
\draw[usual] (0,1) to[out=90,in=270] (0.5,1.5);
\draw[usual,crossline] (0.5,1) to[out=90,in=270] (0,1.5);
\end{tikzpicture}
,\quad
-\bar{3}\colon
\begin{tikzpicture}[anchorbase,scale=1]
\draw[usual] (0.5,0) to[out=90,in=270] (0,0.5);
\draw[usual,crossline] (0,0) to[out=90,in=270] (0.5,0.5);
\draw[usual] (0.5,0.5) to[out=90,in=270] (0,1);
\draw[usual,crossline] (0,0.5) to[out=90,in=270] (0.5,1);
\draw[usual] (0.5,1) to[out=90,in=270] (0,1.5);
\draw[usual,crossline] (0,1) to[out=90,in=270] (0.5,1.5);
\end{tikzpicture}
.
\end{gather*}
The notion of a rational tangle is often defined geometrically, but we will not need the geometric incarnation. The connection between our definition and the geometric is the following lemma.

\begin{Lemma}\label{L:RationalSummands}
Every rational tangle can be obtained as a finite alternating sequence of the form $a_1+_h\bar a_2+_va_3+_h\bar a_4+_v\dots$ or $\bar a_1+_ha_2+_v\bar a_3+_ha_4+_v\dots$ for some $a_i\in\Z$.
\end{Lemma}

\begin{proof}
Well-known (this is essentially Conway's classification of rational tangles), see for example \cite[Section 2.3]{Ad-knots}.
\end{proof}

Let $\mathcal{AL}_{n}$ be the set of {prime} alternating links of $\leq n$ crossings. All probabilities below use this cumulative convention. The cited enumerative results are usually stated for exactly $n$ crossings; the summed form of \autoref{L:EasySandwich} transfers their exponential growth rates to our convention.

\begin{Lemma}\label{L:Alternating1}
We have
\begin{gather*}
\liminf_{n\to\infty}\sqrt[n]{\#\mathcal{AL}_{n}}=
\lim_{n\to\infty}\sqrt[n]{\#\mathcal{AL}_{n}}=\frac{\sqrt{21001}+101}{40}>6.1479.
\end{gather*}
\end{Lemma}

\begin{proof}
The proof is a mild variation of \cite{SuTh-growth}.
Before proving \autoref{L:Alternating1}, let us give an example of the proof strategy in a baby case.

\begin{Example}\label{E:Fib}
The Fibonacci numbers satisfy the recursion $F_0=0$, $F_1=1$ and $F_{n+1}=F_{n}+F_{n-1}$ for $n\geq 1$. Let $F(z)=\sum_{n\geq 0}F_nz^n$ be the generating function. The recursion (essentially immediately) gives the functional equation $(1-z-z^2)F(z)=z$ with poles at $x=1/\phi$ and $y=-\phi$, for the golden ratio $\phi$. Since $x$ is the smallest pole, classical theory then implies $\lim_{n\to\infty}\sqrt[n]{F_n}=x^{-1}=\phi$. Summarized:
\begin{gather*}
\fbox{\mystrut$\text{recursion}\Rightarrow\text{functional equation}
\Rightarrow\text{smallest singularity}\Rightarrow\text{growth rate}$.}
\end{gather*}
(Here and in what follows, smallest for singularities is meant in terms of absolute value.) We apply the same strategy throughout.
\end{Example}

With \autoref{E:Fib} in mind, the proof can be summarized as:
\begin{gather*}
\fbox{\mystrut$\text{planar algebra/operad}\Rightarrow
\text{recursion}\Rightarrow\text{functional equation}
\Rightarrow\text{smallest singularity}\Rightarrow\text{growth rate}$.}
\end{gather*}
To get started, recall from e.g. \cite{SuTh-growth} that:	
\begin{enumerate}

\item Taking the $n$-th root, there is no difference between the growth rate of prime alternating links and prime alternating tangles.

\item By the flyping theorem \cite{MeTh-flyping}, any two reduced alternating diagrams of the same prime alternating link are related by a finite number of operations of the form
\begin{gather*}
\begin{tikzpicture}[anchorbase,scale=1]
\draw[usual] (0.25,0.4) to (0.25,0.6) to[out=90,in=270] (-0.25,1.1);
\draw[usual,tomato,crossline] (-0.25,0.4) to (-0.25,0.6) to[out=90,in=270] (0.25,1.1);
\draw[usual] (0.25,-0.6) to (0.25,0.1);
\draw[usual] (-0.25,-0.6) to (-0.25,0.1);
\draw[fill=cream] (0,0) circle (0.5cm);
\node at (0,0) {$T$};
\end{tikzpicture}
=
\begin{tikzpicture}[anchorbase,scale=1]
\draw[usual,tomato] (0.25,0.4) to (0.25,0.6);
\draw[usual] (-0.25,0.4) to (-0.25,0.6);
\draw[usual] (-0.25,0.1)to (-0.25,-0.6) to[out=270,in=90] (0.25,-1.1);
\draw[usual,crossline] (0.25,0.1) to (0.25,-0.6) to[out=270,in=90] (-0.25,-1.1);
\draw[fill=cream] (0,0) circle (0.5cm);
\node at (0,0) {\reflectbox{$T$}};
\end{tikzpicture},
\end{gather*}
and all variations of this picture. These operations are often called \emph{flypes}.

\end{enumerate}	

For the proof, recall the \emph{cookie-cutter} operation. This operation takes a link diagram and produces a four-valent planar graph with crossings by replacing some crossings with a vertex, e.g.:
\begin{gather*}
\scalebox{0.8}{$\begin{tikzpicture}[anchorbase,scale=1]
\draw[usual] (0.5,0) to[out=90,in=270] (0,0.5);
\draw[usual,crossline] (0,0) to[out=90,in=270] (0.5,0.5);
\draw[usual] (1,0) to[out=90,in=270] (1,0.5);
\draw[usual] (0,0.5) to[out=90,in=270] (0,1);
\draw[usual] (0.5,0.5) to[out=90,in=270] (1,1);
\draw[usual,crossline] (1,0.5) to[out=90,in=270] (0.5,1);
\draw[usual] (0.5,1) to[out=90,in=270] (0,1.5);
\draw[usual,crossline] (0,1) to[out=90,in=270] (0.5,1.5);
\draw[usual] (1,1) to[out=90,in=270] (1,1.5);
\draw[usual] (0,1.5) to[out=90,in=270] (0,2);
\draw[usual] (0.5,1.5) to[out=90,in=270] (1,2);
\draw[usual,crossline] (1,1.5) to[out=90,in=270] (0.5,2);
\draw[usual] (1,2) to[out=90,in=180] (1.25,2.25) to[out=0,in=90] (1.5,2) to (1.5,0) to[out=270,in=0] (1.25,-0.25) to[out=180,in=270] (1,0);
\draw[usual] (0.5,2) to[out=90,in=180] (1.25,2.5) to[out=0,in=90] (2,2) to (2,0) to[out=270,in=0] (1.25,-0.5) to[out=180,in=270] (0.5,0);
\draw[usual] (0,2) to[out=90,in=180] (1.25,2.75) to[out=0,in=90] (2.5,2) to (2.5,0) to[out=270,in=0] (1.25,-0.75) to[out=180,in=270] (0,0);
\end{tikzpicture}
\xrightarrow[\text{crossings}]{\text{forget}}
\begin{tikzpicture}[anchorbase,scale=1]
\draw[usual] (0.5,0) to[out=90,in=270] (0,0.5);
\draw[usual] (0,0) to[out=90,in=270] (0.5,0.5);
\draw[usual] (1,0) to[out=90,in=270] (1,0.5);
\draw[usual] (0,0.5) to[out=90,in=270] (0,1);
\draw[usual] (0.5,0.5) to[out=90,in=270] (1,1);
\draw[usual] (1,0.5) to[out=90,in=270] (0.5,1);
\draw[usual] (0.5,1) to[out=90,in=270] (0,1.5);
\draw[usual] (0,1) to[out=90,in=270] (0.5,1.5);
\draw[usual] (1,1) to[out=90,in=270] (1,1.5);
\draw[usual] (0,1.5) to[out=90,in=270] (0,2);
\draw[usual] (0.5,1.5) to[out=90,in=270] (1,2);
\draw[usual] (1,1.5) to[out=90,in=270] (0.5,2);
\draw[usual] (1,2) to[out=90,in=180] (1.25,2.25) to[out=0,in=90] (1.5,2) to (1.5,0) to[out=270,in=0] (1.25,-0.25) to[out=180,in=270] (1,0);
\draw[usual] (0.5,2) to[out=90,in=180] (1.25,2.5) to[out=0,in=90] (2,2) to (2,0) to[out=270,in=0] (1.25,-0.5) to[out=180,in=270] (0.5,0);
\draw[usual] (0,2) to[out=90,in=180] (1.25,2.75) to[out=0,in=90] (2.5,2) to (2.5,0) to[out=270,in=0] (1.25,-0.75) to[out=180,in=270] (0,0);
\end{tikzpicture}
\xrightarrow[\text{vertex size}]{\text{increase}}
\begin{tikzpicture}[anchorbase,scale=1]
\draw[usual] (0.5,0) to[out=90,in=270] (0,0.5);
\draw[usual] (0,0) to[out=90,in=270] (0.5,0.5);
\draw[usual] (1,0) to[out=90,in=270] (1,0.5);
\draw[usual,fill=cream] (0.25,0.25) circle (0.2cm);
\draw[usual] (0,0.5) to[out=90,in=270] (0,1);
\draw[usual] (0.5,0.5) to[out=90,in=270] (1,1);
\draw[usual] (1,0.5) to[out=90,in=270] (0.5,1);
\draw[usual,fill=cream] (0.75,0.75) circle (0.2cm);
\draw[usual] (0.5,1) to[out=90,in=270] (0,1.5);
\draw[usual] (0,1) to[out=90,in=270] (0.5,1.5);
\draw[usual] (1,1) to[out=90,in=270] (1,1.5);
\draw[usual,fill=cream] (0.25,1.25) circle (0.2cm);
\draw[usual] (0,1.5) to[out=90,in=270] (0,2);
\draw[usual] (0.5,1.5) to[out=90,in=270] (1,2);
\draw[usual] (1,1.5) to[out=90,in=270] (0.5,2);
\draw[usual,fill=cream] (0.75,1.75) circle (0.2cm);
\draw[usual] (1,2) to[out=90,in=180] (1.25,2.25) to[out=0,in=90] (1.5,2) to (1.5,0) to[out=270,in=0] (1.25,-0.25) to[out=180,in=270] (1,0);
\draw[usual] (0.5,2) to[out=90,in=180] (1.25,2.5) to[out=0,in=90] (2,2) to (2,0) to[out=270,in=0] (1.25,-0.5) to[out=180,in=270] (0.5,0);
\draw[usual] (0,2) to[out=90,in=180] (1.25,2.75) to[out=0,in=90] (2.5,2) to (2.5,0) to[out=270,in=0] (1.25,-0.75) to[out=180,in=270] (0,0);
\end{tikzpicture}$}
\,.
\end{gather*}
The result is often called a \emph{template} and the cookie-cutter has produced \emph{slots}, and the same procedure (under the same name) applies to tangles. For both, templates and tangles, we can and will use the sum operation in \autoref{Eq:Sums} for them as well. Moreover, a template is called \emph{rational} if it comes from a rational tangle. The trivial template is the rational template obtained by cookie-cutting the crossing of the tangle $1$. Here is an example of a rational template:
\begin{gather*}
\text{rational, two slots and three crossings: }
\begin{tikzpicture}[anchorbase,scale=1]
\draw[usual] (0.5,0) to[out=90,in=270] (0,0.5);
\draw[usual,crossline] (0,0) to[out=90,in=270] (0.5,0.5);
\draw[usual] (0.5,0.5) to[out=90,in=270] (0,1);
\draw[usual,crossline] (0,0.5) to[out=90,in=270] (0.5,1);
\draw[usual] (0.5,1) to[out=90,in=270] (0,1.5);
\draw[usual,crossline] (0,1) to[out=90,in=270] (0.5,1.5);
\draw[usual] (0,0) to[out=210,in=0] (-0.5,0.5);
\draw[usual] (0,1.5) to[out=150,in=0] (-0.5,1);
\draw[usual] (-0.5,0.5) to[out=180,in=0] (-1,1);
\draw[usual,crossline] (-0.5,1) to[out=180,in=0] (-1,0.5);
\draw[usual] (-1,0.5) to[out=180,in=0] (-1.5,1);
\draw[usual] (-1,1) to[out=180,in=0] (-1.5,0.5);
\draw[fill=cream] (-1.25,0.75) circle (0.2cm);
\draw[fill=cream] (0.25,0.25) circle (0.2cm);
\end{tikzpicture}
.
\end{gather*}
Note that a template can still have crossings. A (tangle) template is \emph{basic} if it has no crossings, more than one slot, is prime in the sense of the operation \autoref{Eq:Sums}, and any circle meeting the template in four points encloses either the whole template or a single slot. Here is an example and a nonexample:
\begin{gather*}
\text{Basic, five slots: }
\begin{tikzpicture}[anchorbase,scale=1]
\draw[usual] (0,0) circle (1cm);
\draw[usual] (-1.4,0) to (1.4,0);
\draw[usual] (0,-1.4) to (0,1.2);
\draw[usual,tomato] (0,1.2) to (0,1.4);
\draw[usual,fill=cream] (0,0) circle (0.2cm);
\draw[usual,fill=cream] (0,1) circle (0.2cm);
\draw[usual,fill=cream] (1,0) circle (0.2cm);
\draw[usual,fill=cream] (0,-1) circle (0.2cm);
\draw[usual,fill=cream] (-1,0) circle (0.2cm);
%%
%	\draw[usual] (1,2) to[out=90,in=180] (1.25,2.25) to[out=0,in=90] (1.5,2) to (1.5,0) to[out=270,in=0] (1.25,-0.25) to[out=180,in=270] (1,0);
%	\draw[usual] (0.5,2) to[out=90,in=180] (1.25,2.5) to[out=0,in=90] (2,2) to (2,0) to[out=270,in=0] (1.25,-0.5) to[out=180,in=270] (0.5,0);
%	\draw[usual] (0,2) to[out=90,in=180] (1.25,2.75) to[out=0,in=90] (2.5,2) to (2.5,0) to[out=270,in=0] (1.25,-0.75) to[out=180,in=270] (0,0);
\end{tikzpicture}
,\quad
\text{not basic, two slots: }
\begin{tikzpicture}[anchorbase,scale=1]
\draw[usual] (-0.9,0.1) to (0.7,0.1);
\draw[usual,tomato] (0.7,0.1) to (0.9,0.1);
\draw[usual] (-0.9,-0.1) to (0.9,-0.1);
\draw[usual,fill=cream] (-0.5,0) circle (0.2cm);
\draw[usual,fill=cream] (0.5,0) circle (0.2cm);
\end{tikzpicture}
.
\end{gather*}
Let $bt_n$ be the number of basic templates with $n$ slots. The following is Tutte's miracle.

\begin{Lemma}\label{L:bt}
The generating function for $bt_n$ is
\begin{gather*}
bt(z)
=
\frac{(1-4z)^{3/2}(z+1)-2z^5-10z^4-10z^3+5z-1}{2(z+1)(z+2)^3}
.
\end{gather*}
The smallest singularity is at $z=1/4$.
\end{Lemma}

\begin{proof}
The key trick is to see that basic templates and so-called rooted c-nets are the same, and the latter have been counted  in \cite{Tu-planar}. From the formula for $bt(z)$ it is then straightforward to see that the singularities are at $z=1/4$, $z=-1$, and $z=-2$.
\end{proof}

Now, the recursion: 

\begin{Lemma}\label{L:Recursion}
Any tangle diagram can be uniquely constructed by repeatedly inserting templates into designated slots, ensuring that any template inserted into a rational template is a basic template.	
\end{Lemma}

\begin{proof}
This works as follows. Call a circle that meets the tangle in four points maximal if it does not enclose the whole tangle or a single crossing, and does not enclose a nontrivial sum of tangles under the operation \autoref{Eq:Sums}. Take the interiors of maximal circles in a tangle diagram, and these tangle diagrams are associated with basic templates, by construction, which can be again obtained by cutting along maximal cycles.
\end{proof}

\begin{Lemma}\label{L:Flype}
Every flype that can occur in a tangle diagram can be traced back to a flype of a rational template somewhere in the hierarchy described in \autoref{L:Recursion} above.
\end{Lemma}

\begin{proof}
See \cite{MeTh-flyping,SuTh-growth}. In the hierarchy from \autoref{L:Recursion}, take the smallest member containing the flype circle. The maximality condition on the basic pieces forces this member to be rational.
\end{proof}

\autoref{L:Flype} and (b) above taken together mean that once we have identified all distinct flype-equivalence classes of rational templates, we can ignore flypes entirely.

Let $at(z)$ be the generating function of the number of alternating prime tangles.
Taking all these observations together gives the functional equation
\begin{gather}\label{Eq:Atrecursion}
at(z)=rt\Big(bt\big(at(z)\big),z\Big),
\end{gather}
where $rt(y,z)$ is the generating function of rational templates up to flype equivalence, $z$ encodes the number of crossings and $y$ the number of slots so that $rt_{m,n}$, the coefficient of $y^mz^n$, is the number of them with $m$ slots and $n$ crossings. Since $bt(z)$ was identified in \autoref{L:bt}, it remains to get a handle on $rt(y,z)$.

\begin{Lemma}\label{L:rt}
The generating function for $rt_{m,n}$ satisfies
\begin{gather*}
2rt(y,z)
=
1+z-y-\Big((1-z+y)^2-\frac{8(z^2-yz+y)}{1-z}\Big)^{1/2}
.
\end{gather*}
\end{Lemma}

\begin{proof}
A bit longish, but nicely explained in \cite{SuTh-growth}.
\end{proof}

One can now use \autoref{L:bt}, \autoref{L:rt} and \autoref{Eq:Atrecursion}, and get the equation
\begin{gather*}
\frac{135 z^2+101 z-20}{108 (z-1)}=0
\;\Rightarrow\;
135 z^2+101 z-20=0,
\end{gather*}
for the singularities of $at(z)$. The statement follows by (a). This finishes the proof of \autoref{L:Alternating1}.
\end{proof}

\begin{Remark}\label{R:Knots}(\emph{Why the restriction to alternating links.})
The proof of \autoref{L:Alternating1} uses the fact that reduced alternating diagrams can be related to alternating links with $n$ crossings. No comparable result is known for arbitrary link diagrams: for a fixed crossing number, there are typically far fewer links than diagrams. Without such a result, the strategy above would give only a $\limsup$ rather than a $\lim$. This is enough for \autoref{L:Alternating2} only because \autoref{L:Alternating1} gives an honest limit.

Similarly, Tutte's miracle \autoref{L:bt} cannot be used to count alternating knots (the number of components is invisible for basic templates), and we do not know any analog of \autoref{L:Alternating1} for knots. The closest result we know is \cite[Theorem 5]{Ch-knots}, which however is not strong enough for our count.

These are the main reasons to restrict to alternating links.
\end{Remark}

Conway mutation along a tangle, \emph{mutation} for short, consists of cutting a tangle from a link and applying one of the following four symmetries of the square. We use the letter $F$ because it is very asymmetric:
\begin{gather*}
\text{id}:
\begin{tikzpicture}[anchorbase,scale=1]
\draw[usual,tomato] (0.25,0.4) to (0.25,0.6);
\draw[usual] (-0.25,0.4) to (-0.25,0.6);
\draw[usual] (0.25,-0.6) to (0.25,0.1);
\draw[usual] (-0.25,-0.6) to (-0.25,0.1);
\draw[fill=cream] (0,0) circle (0.5cm);
\node at (0,0) {$F$};
\end{tikzpicture}
,\quad
\text{reflect left-to-right}:
\begin{tikzpicture}[anchorbase,scale=1]
\draw[usual,tomato] (0.25,0.4) to (0.25,0.6);
\draw[usual] (-0.25,0.4) to (-0.25,0.6);
\draw[usual] (0.25,-0.6) to (0.25,0.1);
\draw[usual] (-0.25,-0.6) to (-0.25,0.1);
\draw[fill=cream] (0,0) circle (0.5cm);
\node at (0,0) {\reflectbox{$F$}};
\end{tikzpicture}
%,\quad
%\begin{tikzpicture}[anchorbase,scale=1]
%\draw[usual,tomato] (0.25,0.4) to (0.25,0.6);
%\draw[usual] (-0.25,0.4) to (-0.25,0.6);
%\draw[usual] (0.25,-0.6) to (0.25,0.1);
%\draw[usual] (-0.25,-0.6) to (-0.25,0.1);
%%%
%\draw[fill=cream] (0,0) circle (0.5cm);
%%%
%\node at (0,0) {\rotatebox{90}{$F$}};
%\end{tikzpicture}
%,\quad
%\begin{tikzpicture}[anchorbase,scale=1]
%\draw[usual,tomato] (0.25,0.4) to (0.25,0.6);
%\draw[usual] (-0.25,0.4) to (-0.25,0.6);
%\draw[usual] (0.25,-0.6) to (0.25,0.1);
%\draw[usual] (-0.25,-0.6) to (-0.25,0.1);
%%%
%\draw[fill=cream] (0,0) circle (0.5cm);
%%%
%\node at (0,0) {\reflectbox{\rotatebox{90}{$F$}}};
%\end{tikzpicture}
,\quad
\text{rotate 180 degrees}:
\begin{tikzpicture}[anchorbase,scale=1]
\draw[usual,tomato] (0.25,0.4) to (0.25,0.6);
\draw[usual] (-0.25,0.4) to (-0.25,0.6);
\draw[usual] (0.25,-0.6) to (0.25,0.1);
\draw[usual] (-0.25,-0.6) to (-0.25,0.1);
\draw[fill=cream] (0,0) circle (0.5cm);
\node at (0,0) {\rotatebox{180}{$F$}};
\end{tikzpicture}
,\quad
\text{reflect top-to-bottom}:
\begin{tikzpicture}[anchorbase,scale=1]
\draw[usual,tomato] (0.25,0.4) to (0.25,0.6);
\draw[usual] (-0.25,0.4) to (-0.25,0.6);
\draw[usual] (0.25,-0.6) to (0.25,0.1);
\draw[usual] (-0.25,-0.6) to (-0.25,0.1);
\draw[fill=cream] (0,0) circle (0.5cm);
\node at (0,0) {\reflectbox{\rotatebox{180}{$F$}}};
\end{tikzpicture}
%,\quad
%\begin{tikzpicture}[anchorbase,scale=1]
%\draw[usual,tomato] (0.25,0.4) to (0.25,0.6);
%\draw[usual] (-0.25,0.4) to (-0.25,0.6);
%\draw[usual] (0.25,-0.6) to (0.25,0.1);
%\draw[usual] (-0.25,-0.6) to (-0.25,0.1);
%%%
%\draw[fill=cream] (0,0) circle (0.5cm);
%%%
%\node at (0,0) {\rotatebox{270}{$F$}};
%\end{tikzpicture}
%,\quad
%\begin{tikzpicture}[anchorbase,scale=1]
%\draw[usual,tomato] (0.25,0.4) to (0.25,0.6);
%\draw[usual] (-0.25,0.4) to (-0.25,0.6);
%\draw[usual] (0.25,-0.6) to (0.25,0.1);
%\draw[usual] (-0.25,-0.6) to (-0.25,0.1);
%%%
%\draw[fill=cream] (0,0) circle (0.5cm);
%%%
%\node at (0,0) {\reflectbox{\rotatebox{270}{$F$}}};
%\end{tikzpicture}
,
\end{gather*}
We regard the last operation as a composition of the other two nontrivial operations, and then reinsert the tangle. In the left-to-right reflection, the crossings must also be mirrored:
\begin{gather*}
\text{original}:
\begin{tikzpicture}[anchorbase,scale=1]
\draw[usual] (0.5,0) to[out=90,in=270] (0,0.5);
\draw[usual,crossline] (0,0) to[out=90,in=270] (0.5,0.5);
\end{tikzpicture}
,\quad
\text{correct reflect left-to-right}:
\reflectbox{\begin{tikzpicture}[anchorbase,scale=1]
\draw[usual] (0,0) to[out=90,in=270] (0.5,0.5);
\draw[usual,crossline] (0.5,0) to[out=90,in=270] (0,0.5);
\end{tikzpicture}}
,\quad
\text{wrong reflect left-to-right}:
\reflectbox{\begin{tikzpicture}[anchorbase,scale=1]
\draw[usual] (0.5,0) to[out=90,in=270] (0,0.5);
\draw[usual,crossline] (0,0) to[out=90,in=270] (0.5,0.5);
\end{tikzpicture}}
.
\end{gather*}
A better, but three-dimensional, way of saying this is (in the same order as above):
\begin{gather*}
\text{id}:
\begin{tikzpicture}[anchorbase,scale=1]
\draw[usual,tomato] (0.25,0.4) to (0.25,0.6);
\draw[usual] (-0.25,0.4) to (-0.25,0.6);
\draw[usual] (0.25,-0.6) to (0.25,0.1);
\draw[usual] (-0.25,-0.6) to (-0.25,0.1);
\draw[fill=cream] (0,0) circle (0.5cm);
\node at (0,0) {$F$};
\end{tikzpicture}
,\quad
\text{z-axis rotate}:
\begin{tikzpicture}[anchorbase,scale=1]
\draw[thick,densely dotted] (0,-0.5) to (0,0.5);
\draw[usual,orchid,directed=0.95] (-0.25,0) to[out=270,in=180] (0,-0.25) to[out=0,in=270] (0.25,0);
\end{tikzpicture}
,\quad
\text{x-axis rotate}:
\begin{tikzpicture}[anchorbase,scale=1]
\draw[thick,densely dotted] (-0.5,0) to (0.5,0);
\draw[usual,orchid,directed=0.95] (0,-0.25) to[out=0,in=270] (0.25,0) to[out=90,in=0] (0,0.25);
\end{tikzpicture}
,\quad
\text{y-axis rotate}:
\begin{tikzpicture}[anchorbase,scale=1]
\node at (0,0) {$\bullet$};
\draw[usual,orchid,directed=0.95] (0.25,0) to[out=90,in=0] (0,0.25) to[out=180,in=90] (-0.25,0);
\end{tikzpicture}
,
\end{gather*}
with standard coordinate system. All of these are \emph{mutants} of the original.

A mutation is \emph{oriented}, or \textit{positive}, if the orientations induced at the four boundary points agree before and after the mutation. For a fixed boundary orientation, two of the four operations above preserve it.

\begin{Example}
A good example of oriented mutation to keep in mind is the following:
\begin{gather*}\label{Eq:Mutation1}
\text{original}:
\begin{tikzpicture}[anchorbase,scale=1]
\draw[usual] (0,0) to (0,0.5);
\draw[usual,crossline] (0,0.5) circle (0.2cm);
\draw[usual,crossline,->] (0,0.5) to (0,1);
\draw[usual,->] (-0.5,0) to (-0.5,1);
\end{tikzpicture}
,\quad
\text{reflect left-to-right}:
\begin{tikzpicture}[anchorbase,scale=1]
\draw[usual] (0,0) to (0,0.5);
\draw[usual,crossline] (0,0.5) circle (0.2cm);
\draw[usual,crossline,->] (0,0.5) to (0,1);
\draw[usual,->] (0.5,0) to (0.5,1);
\end{tikzpicture}
.
\end{gather*}
Another example, the famous pair of the Kinoshita--Terasaka (\texttt{11n\_42}) and the Conway knot (\texttt{11n\_34}) is an example of oriented mutation; here are the corresponding pictures from KnotScape:
\begin{gather*}\label{Eq:Mutation2}
\begin{tikzpicture}[anchorbase]
\node at (0,0) {\includegraphics[height=2cm]{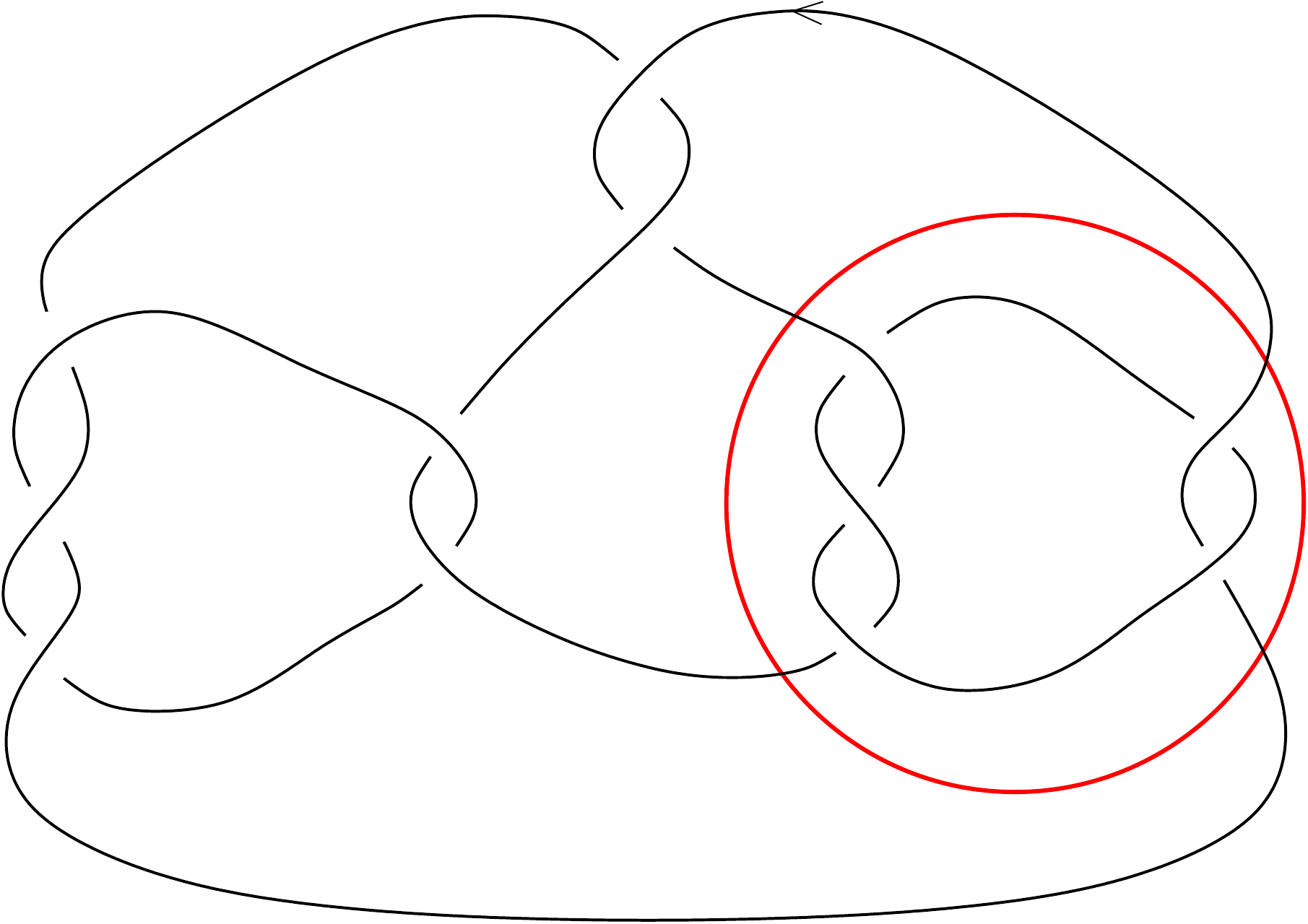}};
\node at (0,-1.2) {11n34};
\end{tikzpicture}
\leftrightarrow
\begin{tikzpicture}[anchorbase]
\node at (0,0) {\includegraphics[height=2cm]{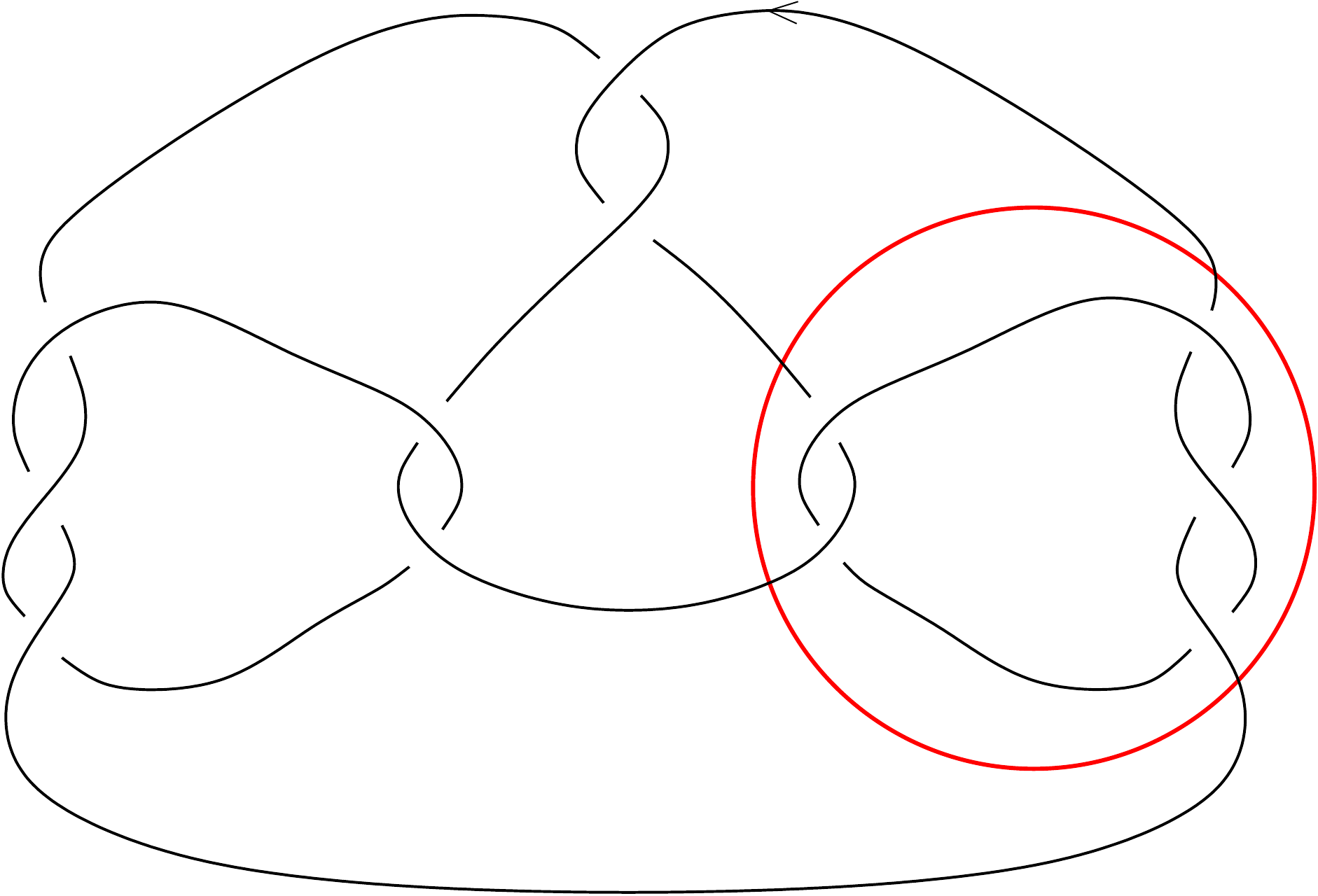}};
\node at (0,-1.2) {11n42};
\end{tikzpicture}
,\quad
\text{mutation}:
\begin{tikzpicture}[anchorbase,scale=1]
\draw[usual,spinach,<-] (0,-.2) to[out=90,in=270] (0.5,0.5);
\draw[usual,crossline,spinach] (0.5,0) to[out=90,in=270] (0,0.5);
\draw[usual,spinach] (0,0.5) to[out=90,in=270] (0.5,1);
\draw[usual,crossline,spinach] (0.5,0.5) to[out=90,in=270] (0,1);
\draw[usual,spinach] (0,1) to[out=90,in=270] (0.5,1.5);
\draw[usual,crossline,spinach] (0.5,1) to[out=90,in=270] (0,1.5) to (0,1.6);	
\draw[usual,tomato] (1.5,0.25) to[out=90,in=270] (1,0.75);
\draw[usual,crossline,spinach] (1,0.25) to[out=90,in=270] (1.5,0.75);
\draw[usual,spinach] (1.5,0.75) to[out=90,in=270] (1,1.25);
\draw[usual,crossline,tomato,->] (1,0.75) to[out=90,in=270] (1.5,1.25) to (1.5,1.35);
\draw[usual,spinach] (0.5,1.5) to[out=90,in=90] (1,1.5) to (1,1.25);
\draw[usual,spinach] (0.5,0) to[out=270,in=270] (1,0) to (1,0.25);
\end{tikzpicture}
\leftrightarrow
\reflectbox{\begin{tikzpicture}[anchorbase,scale=1]
\draw[usual,spinach] (0.5,0) to[out=90,in=270] (0,0.5);
\draw[usual,crossline,spinach] (0,0) to[out=90,in=270] (0.5,0.5);
\draw[usual,spinach] (0.5,0.5) to[out=90,in=270] (0,1);
\draw[usual,crossline,spinach] (0,0.5) to[out=90,in=270] (0.5,1);
\draw[usual,spinach,->] (0.5,1) to[out=90,in=270] (0,1.5) to (0,1.65);
\draw[usual,crossline,spinach] (0,1) to[out=90,in=270] (0.5,1.5);	
\draw[usual,spinach] (1,0.25) to[out=90,in=270] (1.5,0.75);
\draw[usual,crossline,tomato,<-] (1.5,0.05) to[out=90,in=270] (1,0.75);
\draw[usual,tomato] (1,0.75) to[out=90,in=270] (1.5,1.25) to (1.5,1.35);
\draw[usual,crossline,spinach] (1.5,0.75) to[out=90,in=270] (1,1.25);
\draw[usual,spinach] (0.5,1.5) to[out=90,in=90] (1,1.5) to (1,1.25);
\draw[usual,spinach] (0.5,0) to[out=270,in=270] (1,0) to (1,0.25);
\end{tikzpicture}}
.
\end{gather*}
Here the left-hand sides of the knot diagrams are the same, and we mutate along $\bar{3}+_h\bar{2}$ as on the right.
%This is an important difference for several quantum invariants. For example, as shown in \cite{We-mutation}, Khovanov homology can detect the oriented mutation as in \autoref{Eq:Mutation1}, but not as in \autoref{Eq:Mutation2}.
\end{Example}

\begin{Remark}
Not every invariant is insensitive to oriented mutation. Non-rectangular colored HOMFLYPT polynomials can distinguish positive mutants \cite{MR1395780,NaRaSi-colored-mutants}, as can finite type invariants of degree eleven \cite{ChLa-mutants}. Thus the hypothesis in \autoref{T:Decay} is a genuine restriction.
\end{Remark}

Let $\mathcal{ALOM}_{n}$ be the set of prime alternating links of at most $n$ crossings modulo oriented mutation. 

\begin{Lemma}\label{L:Alternating2}
We have
\[
\limsup_{n\to\infty}
\sqrt[n]{\#\mathcal{ALOM}_{n}}
\leq \mu=6.1133171335\ldots<6.114,
\]
where $\mu$ is defined precisely in the proof.
\end{Lemma}

\begin{proof}
The argument is a variation of \autoref{L:Alternating1}, and to make it work we adapt the strategy of \cite{St-number-polynomials}. We start with the following lemma.

\begin{Lemma}\label{L:RationalSummands2}
Oriented mutation can be used to permute the summands (cf. \autoref{L:RationalSummands}) of a rational tangle. The same is true for rational templates.
\end{Lemma}

\begin{proof}
The illustration
\begin{gather*}
-2+_h-\bar 3=
\begin{tikzpicture}[anchorbase,scale=1]
\draw[usual] (0.5,0) to[out=90,in=270] (0,0.5);
\draw[usual,crossline] (0,0) to[out=90,in=270] (0.5,0.5);
\draw[usual] (0.5,0.5) to[out=90,in=270] (0,1);
\draw[usual,crossline] (0,0.5) to[out=90,in=270] (0.5,1);
\draw[usual] (0.5,1) to[out=90,in=270] (0,1.5);
\draw[usual,crossline] (0,1) to[out=90,in=270] (0.5,1.5);
\draw[usual] (0,0) to[out=210,in=0] (-0.5,0.5);
\draw[usual] (0,1.5) to[out=150,in=0] (-0.5,1);
\draw[usual] (-0.5,0.5) to[out=180,in=0] (-1,1);
\draw[usual,crossline] (-0.5,1) to[out=180,in=0] (-1,0.5);
\draw[usual] (-1,0.5) to[out=180,in=0] (-1.5,1);
\draw[usual,crossline] (-1,1) to[out=180,in=0] (-1.5,0.5);
\end{tikzpicture}
\xrightarrow{\text{reflect left-to-right}}
\reflectbox{\begin{tikzpicture}[anchorbase,scale=1]
\draw[usual] (0,0) to[out=90,in=270] (0.5,0.5);
\draw[usual,crossline] (0.5,0) to[out=90,in=270] (0,0.5);
\draw[usual] (0,0.5) to[out=90,in=270] (0.5,1);
\draw[usual,crossline] (0.5,0.5) to[out=90,in=270] (0,1);
\draw[usual] (0,1) to[out=90,in=270] (0.5,1.5);
\draw[usual,crossline] (0.5,1) to[out=90,in=270] (0,1.5);
\draw[usual] (0,0) to[out=210,in=0] (-0.5,0.5);
\draw[usual] (0,1.5) to[out=150,in=0] (-0.5,1);
\draw[usual] (-0.5,1) to[out=180,in=0] (-1,0.5);
\draw[usual,crossline] (-0.5,0.5) to[out=180,in=0] (-1,1);
\draw[usual] (-1,1) to[out=180,in=0] (-1.5,0.5);
\draw[usual,crossline] (-1,0.5) to[out=180,in=0] (-1.5,1);
\end{tikzpicture}}
=
-\bar 3+_h-2,
\end{gather*}
generalizes easily. One can also think of this as rational tangles are determined by a rational number, so this is actually just an isotopy. Now, for our purposes flypes are isotopies, so for rational templates we first use flypes to reorder within $\pm n,\pm \bar{n}$, e.g.:
\begin{gather*}
\begin{tikzpicture}[anchorbase,scale=1]
\draw[usual] (0.5,0) to[out=90,in=270] (0,0.5);
\draw[usual,crossline] (0,0) to[out=90,in=270] (0.5,0.5);
\draw[usual] (0.5,0.5) to[out=90,in=270] (0,1);
\draw[usual,crossline] (0,0.5) to[out=90,in=270] (0.5,1);
\end{tikzpicture}
\xrightarrow[\text{bottom}]{\text{cookie cut}}
\begin{tikzpicture}[anchorbase,scale=1]
\draw[usual] (0.25,0.4) to (0.25,0.6) to[out=90,in=270] (-0.25,1.1);
\draw[usual,crossline] (-0.25,0.4) to (-0.25,0.6) to[out=90,in=270] (0.25,1.1);
\draw[usual] (0.25,-0.6) to (0.25,0.1);
\draw[usual] (-0.25,-0.6) to (-0.25,0.1);
\draw[fill=cream] (0,0) circle (0.5cm);
\node at (0,0) {$T$};
\end{tikzpicture}
=
\begin{tikzpicture}[anchorbase,scale=1]
\draw[usual] (0.25,0.4) to (0.25,0.6);
\draw[usual] (-0.25,0.4) to (-0.25,0.6);
\draw[usual] (-0.25,0.1)to (-0.25,-0.6) to[out=270,in=90] (0.25,-1.1);
\draw[usual,crossline] (0.25,0.1) to (0.25,-0.6) to[out=270,in=90] (-0.25,-1.1);
\draw[fill=cream] (0,0) circle (0.5cm);
\node at (0,0) {\reflectbox{$T$}};
\end{tikzpicture}
\xleftarrow[\text{top}]{\text{cookie cut}}
\begin{tikzpicture}[anchorbase,scale=1]
\draw[usual] (0.5,0) to[out=90,in=270] (0,0.5);
\draw[usual,crossline] (0,0) to[out=90,in=270] (0.5,0.5);
\draw[usual] (0.5,0.5) to[out=90,in=270] (0,1);
\draw[usual,crossline] (0,0.5) to[out=90,in=270] (0.5,1);
\end{tikzpicture}
.
\end{gather*}
If the previous picture contains a slot, say near an extremity, then an actual mutation is needed to pass from left to right:
\begin{gather*}
\begin{tikzpicture}[anchorbase,scale=1]
\draw[usual] (0.5,0) to[out=90,in=270] (0,0.5);
\draw[usual,crossline] (0,0) to[out=90,in=270] (0.5,0.5);
\draw[usual] (0.5,0.5) to[out=90,in=270] (0,1);
\draw[usual,crossline] (0,0.5) to[out=90,in=270] (0.5,1);
\draw[usual] (0.5,1) to[out=90,in=270] (0,1.5);
\draw[usual,crossline] (0,1) to[out=90,in=270] (0.5,1.5);
\draw[usual] (0,0) to[out=210,in=0] (-0.5,0.5);
\draw[usual] (0,1.5) to[out=150,in=0] (-0.5,1);
\draw[usual] (-0.5,0.5) to[out=180,in=0] (-1,1);
\draw[usual,crossline] (-0.5,1) to[out=180,in=0] (-1,0.5);
\draw[usual] (-1,0.5) to[out=180,in=0] (-1.5,1);
\draw[usual] (-1,1) to[out=180,in=0] (-1.5,0.5);
\draw[fill=cream] (-1.25,0.75) circle (0.2cm);
\end{tikzpicture}
\xrightarrow{\text{reflect left-to-right}}
\reflectbox{\begin{tikzpicture}[anchorbase,scale=1]
\draw[usual] (0,0) to[out=90,in=270] (0.5,0.5);
\draw[usual,crossline] (0.5,0) to[out=90,in=270] (0,0.5);
\draw[usual] (0,0.5) to[out=90,in=270] (0.5,1);
\draw[usual,crossline] (0.5,0.5) to[out=90,in=270] (0,1);
\draw[usual] (0,1) to[out=90,in=270] (0.5,1.5);
\draw[usual,crossline] (0.5,1) to[out=90,in=270] (0,1.5);
\draw[usual] (0,0) to[out=210,in=0] (-0.5,0.5);
\draw[usual] (0,1.5) to[out=150,in=0] (-0.5,1);
\draw[usual] (-0.5,1) to[out=180,in=0] (-1,0.5);
\draw[usual,crossline] (-0.5,0.5) to[out=180,in=0] (-1,1);
\draw[usual] (-1,1) to[out=180,in=0] (-1.5,0.5);
\draw[usual] (-1,0.5) to[out=180,in=0] (-1.5,1);
\draw[fill=cream] (-1.25,0.75) circle (0.2cm);
\end{tikzpicture}}
.
\end{gather*}
However this is oriented mutation as the rational tangle type does not change, as we mentioned above. For example,
\begin{gather*}
\begin{tikzpicture}[anchorbase,scale=1]
\draw[fill=cream] (0.25,0.25) circle (0.35cm);
\draw[usual,directed=1] (0.5,0.5) to[out=225,in=335] (0,0.5);
\draw[usual,directed=1] (0,0) to[out=45,in=135] (0.5,0);
\end{tikzpicture}
\text{ works for }
\begin{tikzpicture}[anchorbase,scale=1]
\draw[usual] (0.5,0) to[out=90,in=270] (0,0.5);
\draw[usual,crossline] (0,0) to[out=90,in=270] (0.5,0.5);
\draw[usual] (0.5,0.5) to[out=90,in=270] (0,1);
\draw[usual,crossline] (0,0.5) to[out=90,in=270] (0.5,1);
\draw[usual] (0.5,1) to[out=90,in=270] (0,1.5);
\draw[usual,crossline] (0,1) to[out=90,in=270] (0.5,1.5);
\draw[usual] (0,0) to[out=210,in=0] (-0.5,0.5);
\draw[usual] (0,1.5) to[out=150,in=0] (-0.5,1);
\draw[usual] (-0.5,0.5) to[out=180,in=0] (-1,1);
\draw[usual,crossline] (-0.5,1) to[out=180,in=0] (-1,0.5);
\draw[usual] (-1,0.5) to[out=180,in=0] (-1.5,1);
\draw[usual] (-1,1) to[out=180,in=0] (-1.5,0.5);
\draw[fill=cream] (-1.25,0.75) circle (0.2cm);
\end{tikzpicture}
\text{ and }
\reflectbox{\begin{tikzpicture}[anchorbase,scale=1]
\draw[usual] (0,0) to[out=90,in=270] (0.5,0.5);
\draw[usual,crossline] (0.5,0) to[out=90,in=270] (0,0.5);
\draw[usual] (0,0.5) to[out=90,in=270] (0.5,1);
\draw[usual,crossline] (0.5,0.5) to[out=90,in=270] (0,1);
\draw[usual] (0,1) to[out=90,in=270] (0.5,1.5);
\draw[usual,crossline] (0.5,1) to[out=90,in=270] (0,1.5);
\draw[usual] (0,0) to[out=210,in=0] (-0.5,0.5);
\draw[usual] (0,1.5) to[out=150,in=0] (-0.5,1);
\draw[usual] (-0.5,1) to[out=180,in=0] (-1,0.5);
\draw[usual,crossline] (-0.5,0.5) to[out=180,in=0] (-1,1);
\draw[usual] (-1,1) to[out=180,in=0] (-1.5,0.5);
\draw[usual] (-1,0.5) to[out=180,in=0] (-1.5,1);
\draw[fill=cream] (-1.25,0.75) circle (0.2cm);
\end{tikzpicture}}
.
\end{gather*}
Similarly for the other possibility, and larger cookie cuts.
\end{proof}

By \autoref{L:RationalSummands2}, we may use the full permutation symmetry of the summands at every rational node. This gives a sharper version of Stoimenow's upper bound. (The upper bound from \cite{St-number-polynomials} is $\approx 6.14$, while we get $\approx 6.11$.)

\begin{Lemma}\label{L:rt2}
Put $q(z)=bt(z)$, and let $W(z)$ be the unique formal power series with nonnegative coefficients and zero constant term satisfying
\begin{equation}\label{Eq:MultisetUpperBound}
1+W(z)=\frac{1}{\sqrt{1-z}}
\exp\left(
\frac12\sum_{k\geq1}
\frac{W(z^k)+q\big(W(z^k)\big)}{k}
\right).
\end{equation}
Then the coefficients of $W$ give upper bounds for alternating prime tangles modulo oriented mutation.
\end{Lemma}

\begin{proof}
We use the decomposition from \cite[Proof of Theorems 1 and 2]{St-number-polynomials}, but quotient by every permutation allowed by \autoref{L:RationalSummands2}. Let $\gamma(z)$ count horizontally decomposable rational templates, together with the trivial template, after the slots have been filled, and let $W(z)$ count all resulting tangles. A summand is either a single crossing or an object counted by $\gamma$. Since adjacent transpositions are realized by oriented mutations, the summands form an unlabeled multiset. Consequently,
\[
\gamma=q(W)+\operatorname{MSET}_{\geq2}(z+\gamma),
\qquad
W=2\gamma+z-q(W).
\]
Here
\[
\operatorname{MSET}(A)(z)
=\exp\left(\sum_{k\geq1}\frac{A(z^k)}{k}\right),
\qquad
\operatorname{MSET}_{\geq2}(A)=\operatorname{MSET}(A)-1-A.
\]
The second equation accounts for horizontal and vertical sums and subtracts the twice-counted trivial template. Eliminating $\gamma$ gives \autoref{Eq:MultisetUpperBound}. Every alternating prime tangle has a representative in this construction, while we do not quotient by any additional symmetry at a basic node. Thus the coefficients are upper bounds, as claimed.
\end{proof}

It remains to determine the radius of convergence of $W$. The first singularity of $q$ is $1/4$, by \autoref{L:bt}. Let $r_0$ be the first positive number for which $W(r_0)=1/4$, and set $\mu=r_0^{-1}$. The terms in \autoref{Eq:MultisetUpperBound} with $k\geq2$ are analytic at $r_0$. Moreover,
\[
\frac{1}{1+W}-\frac{1+q'(W)}{2}
\geq
\frac45-\frac12\left(1+\frac{167}{2025}\right)>0
\qquad (0\leq W\leq1/4),
\]
where $q'(1/4)=167/2025$. Hence there is no earlier implicit-function singularity, and $r_0$ is the dominant singularity. Recursive numerical evaluation at $z,z^2,z^3,\ldots$ gives
\[
r_0=0.1635773145998734\ldots,
\qquad
\mu=r_0^{-1}=6.113317133528573\ldots.
\]
For the stated rigorous bound, evaluate the logarithm of \autoref{Eq:MultisetUpperBound} at $z=1/6.114$ and $W=1/4$. Truncating after $k=20$ leaves a positive margin greater than $1.2\cdot10^{-5}$, while the remaining nonnegative tail is less than $3\cdot10^{-18}$. Thus $W(1/6.114)<1/4$, and so $\mu<6.114$. As in \autoref{L:Alternating1}, passing from rooted tangles to links does not increase the exponential growth rate.
\end{proof}

Given an invariant of links $Q$, we define its value ratio on $\mathcal{AL}_n$ by
\begin{gather*}
v^{AL}_Q(n)=
\#\{Q(L) \mid L\in\mathcal{AL}_{n}\}/
\#\mathcal{AL}_{n}.
\end{gather*}
That is, $v^{AL}_Q(n)$ measures the number of distinct values attained by $Q$ on prime alternating links with at most $n$ crossings, divided by the number of such links. We say that $Q$ is \emph{insensitive to oriented mutation} if $Q(L)=Q(L^{\prime})$ whenever $L$ and $L^{\prime}$ are related by a finite sequence of oriented mutations. We also write
\[
m_Q(L,n)=\#\{L'\in\mathcal{AL}_n\mid Q(L')=Q(L)\}
\]
for the size of the $Q$-fiber containing $L$.

\begin{Theorem}\label{T:Decay}
\textup{(}\textbf{Exponential decay.}\textup{)}
Set
\[
\lambda=\frac{\sqrt{21001}+101}{40},\qquad
\mu=r_0^{-1}=6.1133171335\ldots,\qquad
\rho=\frac{\mu}{\lambda}=0.9943699249\ldots,
\]
where $r_0$ is defined by \autoref{Eq:MultisetUpperBound} and $W(r_0)=1/4$.
For any link invariant $Q$ that is insensitive to oriented mutation and every $\delta\in(\rho,1)$, we have
\begin{gather*}
v^{AL}_Q(n)\in O(\delta^{n}).
\end{gather*}
Consequently, the detection ratio also decays exponentially.
\end{Theorem}

\begin{proof}
Since $Q$ is constant on oriented mutation classes,
\[
\#Q(\mathcal{AL}_n)\leq \#\mathcal{ALOM}_n.
\]
Now combine \autoref{L:EasySandwich}, \autoref{L:Alternating1} and \autoref{L:Alternating2}. Their growth constants are $\lambda$ and at most $\mu$, respectively, so every $\delta>\mu/\lambda$ gives the claimed bound. Finally, the links detected by $Q$ form singleton fibers, and hence their number is at most $\#Q(\mathcal{AL}_n)$. Thus the detection ratio is bounded by $v^{AL}_Q(n)$.
\end{proof}

\begin{Theorem}\label{C:Largefibers}
Let $Q$ be insensitive to oriented mutation. For every $\eta\in(1,\rho^{-1})$ there is some $\theta\in(0,1)$ such that
\[
\frac{\#\{L\in\mathcal{AL}_n\mid m_Q(L,n)<\eta^n\}}{\#\mathcal{AL}_n}
\in O(\theta^n).
\]
In particular, with probability tending to one exponentially fast, a prime alternating link shares its $Q$-value with exponentially many other prime alternating links.
\end{Theorem}

\begin{proof}
There are at most $\#Q(\mathcal{AL}_n)$ fibers, and the fibers counted in the numerator have fewer than $\eta^n$ elements. Hence the numerator is at most $\eta^n\#Q(\mathcal{AL}_n)$. Apply \autoref{T:Decay} and choose $\delta\in(\rho,\eta^{-1})$; then one may take $\theta=\eta\delta<1$.
\end{proof}

We need two definitions and a lemma, which are essentially borrowed from \cite{TuZh-knot-data}. We use quantum invariants to mean the same as in that paper (both, polynomial and homological).
\begin{enumerate}[label=(\roman*)]

\item A \emph{skein relation} is a relation of the following form. Let $a,b,c$ be elements in a ring, then
\begin{gather*}
a\cdot
\begin{tikzpicture}[anchorbase,scale=1]
\draw[usual,directed=1] (0.5,0) to[out=90,in=270] (0,0.5);
\draw[usual,crossline,directed=1] (0,0) to[out=90,in=270] (0.5,0.5);
\end{tikzpicture}
+b\cdot
\begin{tikzpicture}[anchorbase,scale=1]
\draw[usual,directed=1] (0,0) to[out=90,in=270] (0.5,0.5);
\draw[usual,crossline,directed=1] (0.5,0) to[out=90,in=270] (0,0.5);
\end{tikzpicture}
+c\cdot
\begin{tikzpicture}[anchorbase,scale=1]
\draw[usual,directed=1] (0,0) to[out=90,in=270] (0,0.5);
\draw[usual,directed=1] (0.5,0) to[out=90,in=270] (0.5,0.5);
\end{tikzpicture}
=
0
,\quad \text{for } a,b\text{ invertible}.
\end{gather*}

\item \emph{Multiplicity freeness} of a quantum invariant is the property that $V\otimes V$ is multiplicity free, where $V$ is the representation used to color the link components.

\end{enumerate}

The following is well-known, but we sketch a proof for completeness.

\begin{Lemma}\label{L:PolyDone}
Any polynomial quantum invariant that satisfies a skein relation or is multiplicity free is invariant under (general) mutation.
\end{Lemma}

\begin{proof}
Given a skein relation, one can proceed as follows. The key idea is that the associated invariant is constant on skein equivalence classes, where two diagrams are skein equivalent if they are related by a sequence of skein moves and isotopies. One can then show that mutation preserves skein equivalence by induction on the number of crossings in the mutated tangle: smoothing a crossing yields a skein triple whose mutated versions remain skein equivalent, allowing the inductive step to go through.

For the multiplicity-free criterion, cut along the Conway sphere and decompose the four-point state space into fusion channels. Each multiplicity space is one-dimensional, so mutation only transposes scalar blocks. Closing with the complementary tangle leaves the trace unchanged. See \cite{MR1395780,MR3381330}, which also explain what can go wrong when multiplicities occur.
\end{proof}

\begin{proof}[Proof of \autoref{T:Main}]
%We start with a standard fact.
%
%\begin{Lemma}\label{L:RationNumbers}
%For any rational tangle $T$ there are rational tangles $U,V$ such that $T+_vU=0$ and $T+_hV=0$.
%\end{Lemma}	
%
%\begin{proof}
%To see this recall (see e.g. \cite[Section 2.3]{Ad-knots}) that rational tangles correspond to and behave like rational numbers (with the number $0$ being unknotted), so we can just compose a rational tangle for $p/q$ with the one for $-p/q$.
%\end{proof}
All of the invariants marked green in \autoref{Tab:Main} are invariant under oriented mutation:
\begin{enumerate}[label=$\blacktriangleright$]

\item For the polynomial invariants, this follows from the skein or multiplicity-free criterion in \autoref{L:PolyDone}. This includes the symmetric and exterior colors listed in the table. The stable HOMFLYPT formulation can also be obtained from \cite[Section 4]{TuVaWe-super-howe}. We make no claim here for non-rectangular colors with multiplicities.

\item For Khovanov homology over $\Ff_2$ or odd Khovanov homology this follows from \cite{We-mutation-2} and \cite{Bl-odd}.

\item The signature, determinant and double branched cover are also invariant under oriented mutation. See \cite{MR3381330} for a summary and \cite{MR3049933} for the branched cover and its lattice. The same is then true of the Heegaard Floer homology of the double branched cover.

%\item For the remaining homological invariants we use \cite[Theorems 1.3 and 1.6]{LaCo-mutation} via \autoref{L:RationNumbers}.

\end{enumerate}
Additionally:
\begin{enumerate}[label=$\blacktriangleright$]

\item For integral Khovanov homology, this follows from the fact that the Khovanov homology of a non-split alternating link is determined by its Jones polynomial and signature (both are invariant under oriented mutation), and a combinatorial gadget that is oriented mutation invariant by \autoref{L:Formula}. This is first proven in \cite[Corollary 1.J]{Sh-khovanov-thin}, where the dependence on the combinatorial formula is implicit, and more explicitly written in \cite[Lemma 3.5]{KeSch-on-comp-comp-kh}.

\item For $\widehat{\mathrm{HFK}}$ over $\Ff_2$, \cite[Theorem 2]{MO08} shows that an alternating link is Floer homologically thin. Hence its bigraded $\widehat{\mathrm{HFK}}$ is determined by the Alexander polynomial and signature, which are invariant under oriented mutation.

\end{enumerate}

\begin{Lemma}\label{L:Formula}
Let $L$ and $L'$ be oriented links related by an oriented mutation.
Then the combinatorial expression for 
$\dim_{\mathbb{Q}} H^{*}_{\mathrm{Lee}}(\placeholder;\mathbb{Q})$
given in \cite[Proposition~4.3]{Lee} takes the same values for $L$
and $L'$.
\end{Lemma}

\begin{proof}
Let us start by dividing the components of $L$ into three disjoint sets $C_{\operatorname{in}}$, $C_{\operatorname{out}}$ and $\{a,b\}$. The set $C_{\operatorname{in}}$ contains the components of $L$ which are completely inside the mutation disk, $C_{\operatorname{out}}$ those that are outside of the disk, and  $\{a,b\}$ has the remaining one or two components. The components of $L$ are identified with the corresponding components of $L'$ with $a$ and $b$ being identified with the respective components of $L'$ which agree with them outside the mutation disk. We denote $\operatorname{lk}_{x,y}$ as the linking number of $x$ and $y$ in $L$ and  $\operatorname{lk}'_{x,y}$ as the linking number of the corresponding components in $L'$. 

The combinatorial expression for the rational Lee homology is given by
$$
\dim_{\mathbb Q}H^{i}_{\mathrm{Lee}}(L;\mathbb{Q})=\# \bigg\{ E\subset \{ \text{components of }L \}\bigg\vert \, \sum_{x\in E,\ y\notin E} 2\operatorname{lk}_{x,y} =i\bigg\}.
$$
To show that this is invariant under oriented mutation, we construct a bijection 
$$
\varphi\colon \big\{ E\subset \{ \text{components of }L \} \big\} \to \big\{ E\subset \{ \text{components of }L' \} \big\}
$$
by
$$
\varphi(E)= \begin{cases}
			E & \text{if } a,b\in E \text{ or } a,b\notin E,\\
            (E\cap (\{a,b\}\cup C_{\operatorname{out}}))\cup (E^c \cap C_{\operatorname{in}}) & \text{otherwise}.
		 \end{cases}   
$$
The linking numbers of $L$ and $L'$ either agree, $\operatorname{lk}_{x,y}=\operatorname{lk}'_{x,y}$, or admit a symmetry $\operatorname{lk}_{a,x}=\operatorname{lk}'_{b,x}$ and $\operatorname{lk}_{b,x}=\operatorname{lk}'_{a,x}$ whenever $x\in  C_{\operatorname{in}}$. 
It follows that 
$$
\sum_{x\in E,\ y\notin E} \operatorname{lk}_{x,y}=\sum_{x\in \varphi (E),\ y\notin \varphi(E)} \operatorname{lk}'_{x,y}
$$
for all $E$ and thus $\dim_{\mathbb Q}H^{i}_{\mathrm{Lee}}(L;\mathbb{Q})=\dim_{\mathbb Q}H^{i}_{\mathrm{Lee}}(L' ;\mathbb{Q})$.
\end{proof}

Now apply \autoref{T:Decay} and \autoref{C:Largefibers}. This proves all three assertions in \autoref{T:Main}.
\end{proof}

Informally, one might summarize the argument as follows:
\begin{gather}\label{Eq:Slogan}
\fcolorbox{orchid!50}{spinach!10}{\mystrut Any local relation on the set of knots forces exponential decay.}
\end{gather}
However, this is not a statement we can prove, so it should be taken with a pinch of salt.
In any case, the TQFT-like fashion in which quantum invariants are defined makes them exciting on the one hand (e.g. for computations, provable statement etc.) but is a fundamental drawback for them being knot invariants on the other hand.

\begin{Remark}\label{R:KR}
With \autoref{Eq:Slogan} in mind, the above strategy should also apply to Khovanov--Rozansky, and HOMFLYPT homology. Integral Khovanov homology is delicate under mutation in general; see \cite{We-mutation,LaCo-mutation,KoWaZi}. For non-split alternating links, the result of Shumakovitch used above reduces it to the Jones polynomial and signature. For HOMFLYPT homology (and partially for Khovanov--Rozansky homology), \cite[Corollary 1.2]{Ja-mutation-homology} establishes insensitivity to oriented mutation in the case of knots. Nevertheless, numerous subtleties remain, some potentially quite difficult, such as restricting results from links to knots (cf. \autoref{R:Knots}), and we therefore chose not to pursue this direction in the present work.
\end{Remark}

\section{Further consequences and sharper forms}\label{S:Consequences}

The proof above gives more than the failure of one invariant at a time. We record four consequences. The first two make the topological content explicit, the third allows controlled mutation sensitivity, and the fourth keeps the polynomial factors supplied by the generating function singularities.

\subsection{Exponentially many mutants}

Write $[L]_{\mathrm{om}}$ for the oriented mutation class of $L$. Taking the oriented mutation class itself as an invariant gives the following.

\begin{Proposition}\label{C:ManyMutants}
\textup{(}\textbf{Exponentially many mutants.}\textup{)}
For every $\eta\in(1,\rho^{-1})$ there is a $\theta\in(0,1)$ such that
\[
\frac{\#\{L\in\mathcal{AL}_n\mid
\#([L]_{\mathrm{om}}\cap\mathcal{AL}_n)<\eta^n\}}
{\#\mathcal{AL}_n}
\in O(\theta^n).
\]
In particular, with probability tending to one exponentially fast, a prime alternating link with at most $n$ crossings has at least $\eta^n$ pairwise inequivalent oriented mutants. One may take any
\[
1<\eta<\rho^{-1}=1.0056619523\ldots.
\]
\end{Proposition}

\begin{proof}
Apply \autoref{C:Largefibers} to the invariant $L\mapsto[L]_{\mathrm{om}}$. Its fibers are precisely the oriented mutation classes, so we are done.
\end{proof}

\subsection{All invariants jointly}

Let $(Q_i)_{i\in I}$ be any family of invariants that are insensitive to oriented mutation, where the index set $I$ need not be finite, and set
\[
\mathbf Q(L)=\big(Q_i(L)\big)_{i\in I}.
\]
The joint invariant $\mathbf Q$ is still constant on every oriented mutation class.

\begin{Proposition}\label{C:JointInvariants}
\textup{(}\textbf{Joint failure of detection.}\textup{)}
The conclusions of \autoref{T:Decay}, \autoref{C:Largefibers} and \autoref{P:PuiseuxRate} hold for the joint invariant $\mathbf Q$. In particular, they hold when all the green invariants in \autoref{Tab:Main} are used simultaneously. Thus even their complete joint package identifies only an exponentially vanishing proportion of prime alternating links.
\end{Proposition}

\begin{proof}
Since $\mathbf Q$ is constant on oriented mutation classes,
\[
\#\mathbf Q(\mathcal{AL}_n)\leq\#\mathcal{ALOM}_n.
\]
The proofs of \autoref{T:Decay} and \autoref{C:Largefibers} now apply without change.
\end{proof}

\subsection{Mild mutation sensitivity}

Strict mutation invariance is not necessary. For an arbitrary invariant $Q$ (not necessarily mutation invariant), define
\[
s_Q(n)=
\max_{C\in\mathcal{ALOM}_n}
\#\{Q(L)\mid L\in C\},
\qquad
\sigma_Q=\limsup_{n\to\infty}s_Q(n)^{1/n}.
\]
Thus $s_Q(n)$ measures how many pieces $Q$ can split one oriented mutation class into.

\begin{Proposition}\label{T:MildMutation}
\textup{(}\textbf{Controlled mutation sensitivity.}\textup{)}
Suppose that $\sigma_Q\rho<1$. Then, for every $\delta\in(\sigma_Q\rho,1)$,
\[
v^{AL}_Q(n),\ d^{AL}_Q(n)\in O(\delta^n).
\]
Moreover, for every $\eta\in(1,(\sigma_Q\rho)^{-1})$, the probability that a uniformly random $L\in\mathcal{AL}_n$ belongs to a $Q$-fiber of size smaller than $\eta^n$ tends to zero exponentially. In particular, the original conclusions remain valid whenever $Q$ takes only subexponentially many values on each oriented mutation class.
\end{Proposition}

\begin{proof}
Every oriented mutation class contributes at most $s_Q(n)$ values of $Q$, so
\[
\#Q(\mathcal{AL}_n)
\leq s_Q(n)\#\mathcal{ALOM}_n.
\]
Together with \autoref{L:Alternating1} and \autoref{L:Alternating2}, this gives
\[
\limsup_{n\to\infty}
\sqrt[n]{v^{AL}_Q(n)}
\leq \sigma_Q\frac{\mu}{\lambda}
=\sigma_Q\rho<1.
\]
The detection ratio is at most the value ratio. Finally, the number of links in $Q$-fibers of size smaller than $\eta^n$ is at most
\[
\eta^n\#Q(\mathcal{AL}_n),
\]
and the last statement follows because $\eta\sigma_Q\rho<1$.
\end{proof}

\subsection{Puiseux asymptotics}

The proof of exponential decay used only the radii of convergence. The series in \autoref{L:Alternating1} is algebraic. The sharper upper-bound series $W$ in \autoref{L:Alternating2} is not algebraic, but \autoref{Eq:MultisetUpperBound} is a subcritical P\'olya equation: the terms with $k\geq2$ are analytic at $r_0$, and the singularity comes from the $3/2$-power in $q(W)$. Thus both counting series have Puiseux expansions at their dominant singularities. Keeping the resulting polynomial factors gives a sharper universal rate.

\begin{Proposition}\label{P:PuiseuxRate}
\textup{(}\textbf{The boundary exponential rate.}\textup{)}
There is a constant $\kappa\geq0$, depending only on the two counting series and not on $Q$, such that every invariant $Q$ insensitive to oriented mutation satisfies
\[
v^{AL}_Q(n),\ d^{AL}_Q(n)
\in O(n^\kappa\rho^n).
\]
For every $\eta\in(1,\rho^{-1})$ we also have
\[
\frac{\#\{L\in\mathcal{AL}_n\mid m_Q(L,n)<\eta^n\}}
{\#\mathcal{AL}_n}
\in O\bigl(n^\kappa(\eta\rho)^n\bigr).
\]
\end{Proposition}

\begin{proof}
Let $U_n$ denote the cumulative coefficient sequence of the upper-bound series $W$ constructed in the proof of \autoref{L:Alternating2}. Thus
\[
\#\mathcal{ALOM}_n\leq U_n.
\]
Newton--Puiseux expansion for the lower series, together with the singular implicit-function argument above and the standard transfer results \cite[Chapter~VI]{FlSe-analytic-combinatorics}, gives constants $c,C>0$ and $a,b\geq0$ such that, for all sufficiently large $n$,
\[
\#\mathcal{AL}_n\geq c n^{-a}\lambda^n,
\qquad
U_n\leq C n^b\mu^n.
\]
Here passing between the rooted tangle series used in the count and links, and then taking cumulative counts, changes only polynomial factors. Hence
\[
v^{AL}_Q(n)
\leq\frac{U_n}{\#\mathcal{AL}_n}
\leq\frac{C}{c}n^{a+b}\rho^n.
\]
Take $\kappa=a+b$. The detection estimate follows from $d^{AL}_Q(n)\leq v^{AL}_Q(n)$. For the last claim, the numerator is at most $\eta^n\#Q(\mathcal{AL}_n)$, and the same estimate gives the stated bound.
\end{proof}

\section{Big data and categorification}\label{S:BigData}

We now closely follow \cite{TuZh-knot-data}, but restrict our attention to the value and detection ratios for knots.

\begin{Remark}\label{R:Code}
All data files and additional material
(such as interactive playgrounds, higher resolution pictures, code, a possibly empty erratum and much more) can be found online at \cite{LaTuVaZh-quantum-big-data-code}. The online material separates alternating and non-alternating knots, includes additional invariants and, in some cases, extends to 19 crossings.
In addition, several of the large-scale analyses of knot invariants have been delegated to these repositories. 
This approach comes with both advantages and drawbacks, but the repositories are intended to remain dynamic and evolve over time, which motivated our choice. We strongly encourage readers to consult the companion website
\begin{center}
\url{https://dustbringer.github.io/web--knot-invariant-comparison/}.
\end{center}
It contains the full histograms of fiber sizes, including singleton fibers, as well as successive quotients and filters for alternating knots, non-alternating knots and mutation classes. It also contains data for many other invariants, plots of computation time, rank and roots, and exploratory views using ball mapper, persistent homology and nearest neighbors. The five plots below are only a summary.
\end{Remark}

The data discussion has three steps. We first describe the invariants, computations and available crossing ranges. We then compare the value and detection ratios. Finally, we use the Jones distribution to show what these two summary statistics leave out.

The invariants we consider are the following (the homologies always categorify the associated polynomials). The reader unfamiliar with these can find more details in the references explaining the computation used. The databases used are \cite{knotinfo} for 3-16 crossings and \cite{knots} for 17-18 crossings.
\begin{enumerate}

\item Khovanov homology and odd Khovanov homology, both categorifying the Jones polynomial. We abbreviate these as K, O, and J, respectively.
\begin{enumerate}

\item K was computed using the \texttt{KnotTheory} Mathematica package from \cite{KnotTheory}.

\item O was computed using the \texttt{KnotJob} program from \cite{KnotJob}.

\item J was computed using a homemade program available at \cite{LaTuVaZh-quantum-big-data-code}. We used it instead of existing calculators, such as Regina, because it was faster on our database.

\end{enumerate}
\vspace{0.3em}
For K and O we computed the homologies with rational coefficients. Their Poincar{\'e} polynomials lie in $\N[q,t]$, and specialization at $t=-1$ gives J. We also consider the specialization $t=1$ of K, called KT1. The data cover knots with $3$--$18$ crossings.

\item Knot Floer (HFK) homology and the Alexander polynomial, denoted HFK2 and A.
\begin{enumerate}

\item HFK was computed using the Knot Floer homology calculator, \cite{Knot-Floer-homology-calculator}, with algorithms from \cite{MR4002230} and \cite{szabo2019algebrasmatchingsknotfloer}.

\item A was computed using the \texttt{KnotTheory} Mathematica package from \cite{KnotTheory}.

\end{enumerate}
\vspace{0.3em}
HFK was computed over $\Ff_2$, giving a polynomial in $\N[q,t]$. We also consider its $t=1$ specialization; the $t=-1$ specialization is A. The data cover knots with $3$--$18$ crossings.

\item Khovanov--Rozansky homology for $n=3$, categorifying the $\mathfrak{sl}_3$ polynomial. We denote these by 
KR3 (homology) and SL3 (polynomial).
\begin{enumerate}

\item KR3 was computed using the \texttt{KnotJob} program from \cite{KnotJob}.

\item SL3 was computed using the \texttt{KnotTheory} Mathematica package from \cite{KnotTheory}.

\end{enumerate}
\vspace{0.3em}
We use the same conventions as for K. The KR3 data cover knots with $3$--$15$ crossings and the SL3 data cover knots with $3$--$18$ crossings.

\item HOMFLYPT homology and its decategorification, the HOMFLYPT polynomial. These are denoted HH (homology) and H (polynomial).
\begin{enumerate}

\item HH was computed with the program used in \cite{homfly}, available via the link in \cite{homfly}. In contrast to the other invariants, this was computed using the knot database in \cite{Knot-database} as it provides efficient braid presentations of the knots (which is very important for the used program). This is the bottleneck of the calculation: it ran for $\approx$215 days on the server of 
the Laboratoire de Math{\'e}matiques Blaise Pascal.

\item H was computed using the \texttt{KnotTheory} Mathematica package, \cite{KnotTheory}.

\end{enumerate}
We use the same conventions as for K. The HH data cover knots with $3$--$13$ crossings and the H data cover knots with $3$--$16$ crossings. Moreover, the data set is incomplete, as detailed in \autoref{table:percent-uniq-val-HOMFLY}.

\end{enumerate}

\begin{Remark}
To the best of our knowledge, these homologies have not been computed on this scale before.
\end{Remark}

We are interested in how many distinct values these invariants take on the list of prime knots $\mathcal{K}_{n}$ (with the same conventions as in \cite{TuZh-knot-data}, e.g. not including mirrors and for $\leq n$ crossings), and we measure this as a ratio:
\begin{gather*}
v^K_Q(n)=
\#\{Q(K) \mid K\in\mathcal{K}_{n}\}/
\#\mathcal{K}_{n}.
\end{gather*}
This is the value ratio, not the detection ratio. The latter is
\[
d^K_Q(n)=
\frac{\#\big\{K\in\mathcal K_n\mid Q^{-1}\big(Q(K)\big)\cap\mathcal K_n=\{K\}\big\}}
{\#\mathcal K_n},
\qquad d^K_Q(n)\leq v^K_Q(n).
\]
The distinction matters: $v^K_Q(n)$ counts fibers, whereas $d^K_Q(n)$ counts knots lying in singleton fibers.

\begin{figure}[!ht]
\begin{center}
    % zoom is 150% on firefox
    \includegraphics[width=0.6\textwidth]{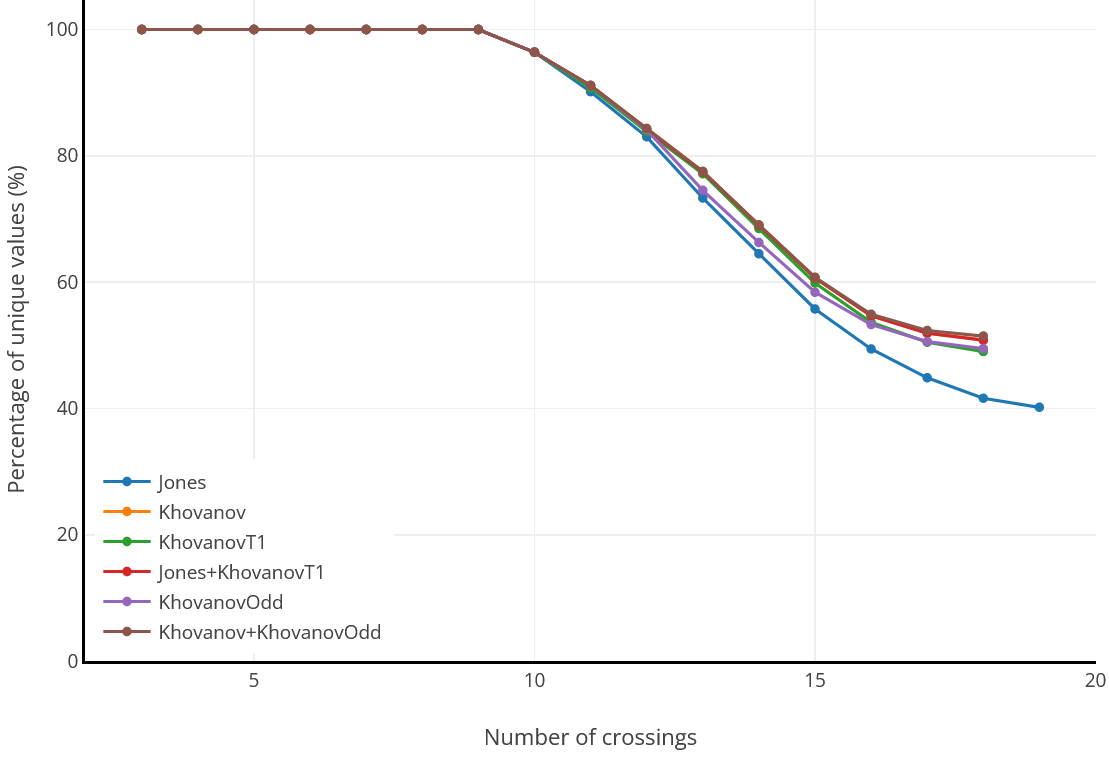}
\end{center}
\caption{Value ratios for the Jones family (equivalently, the number of distinct values as a percentage of the number of knots).}
\label{fig:percent-uniq-val-K}
\end{figure}

\begin{figure}[!ht]
\begin{center}
    \includegraphics[width=0.6\textwidth]{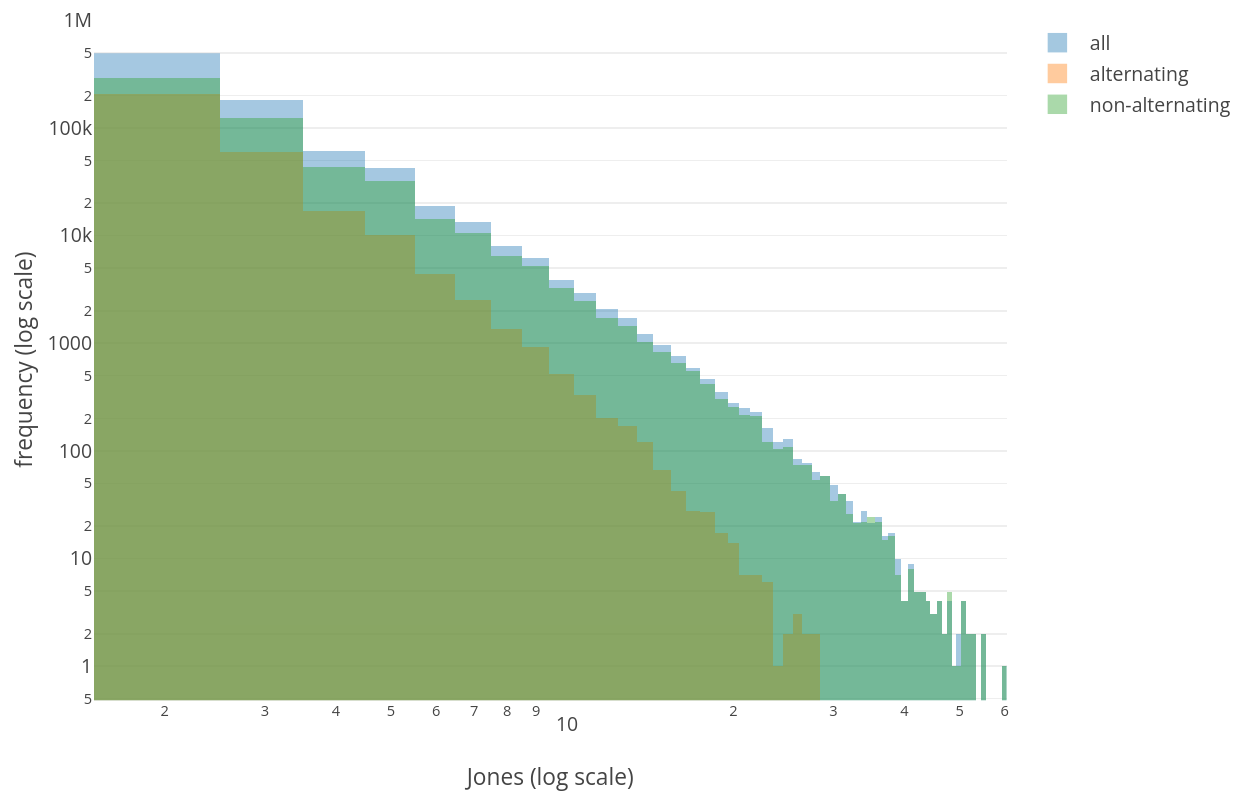}
\end{center}
\caption{Log-log distribution of Jones fiber sizes for prime knots through 16 crossings, with mirrors not distinguished. The horizontal axis records fiber size and the vertical axis the number of fibers of that size. The nearly linear profiles for the full, alternating and non-alternating data are consistent with an approximately power-law tail.}
\label{fig:jones-distribution}
\end{figure}

\begin{figure}[!ht]
\begin{center}
    \includegraphics[width=0.6\textwidth]{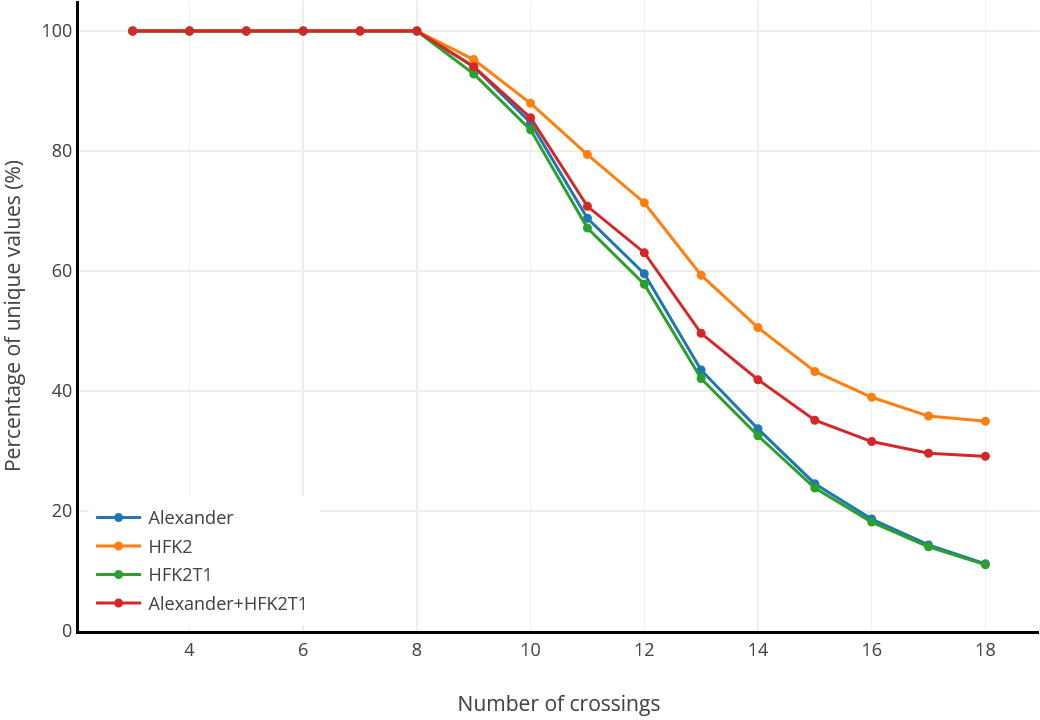}
\end{center}
\caption{Value ratios for the Alexander family.}
\label{fig:percent-uniq-val-HFK}
\end{figure}

\begin{figure}[!ht]
\begin{center}
    \includegraphics[width=0.6\textwidth]{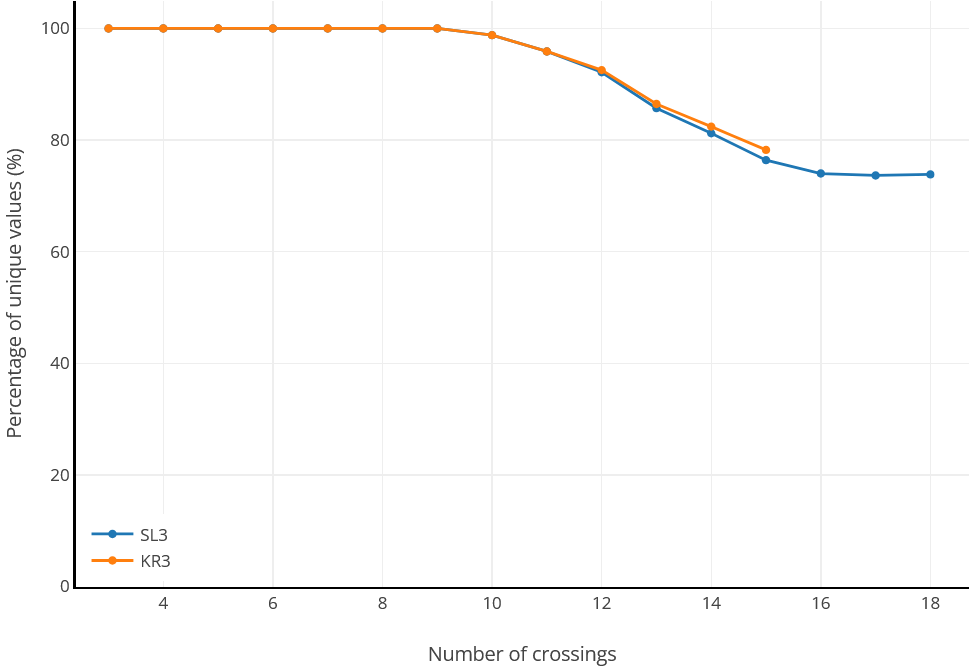}
\end{center}
\caption{Value ratios for the SL3 family. The comparison of specializations is not illustrated as the homology is essentially its Euler characteristic polynomial on the data set.}
\label{fig:percent-uniq-val-SL3}
\end{figure}

\begin{figure}[!ht]
\begin{center}
    \includegraphics[width=0.6\textwidth]{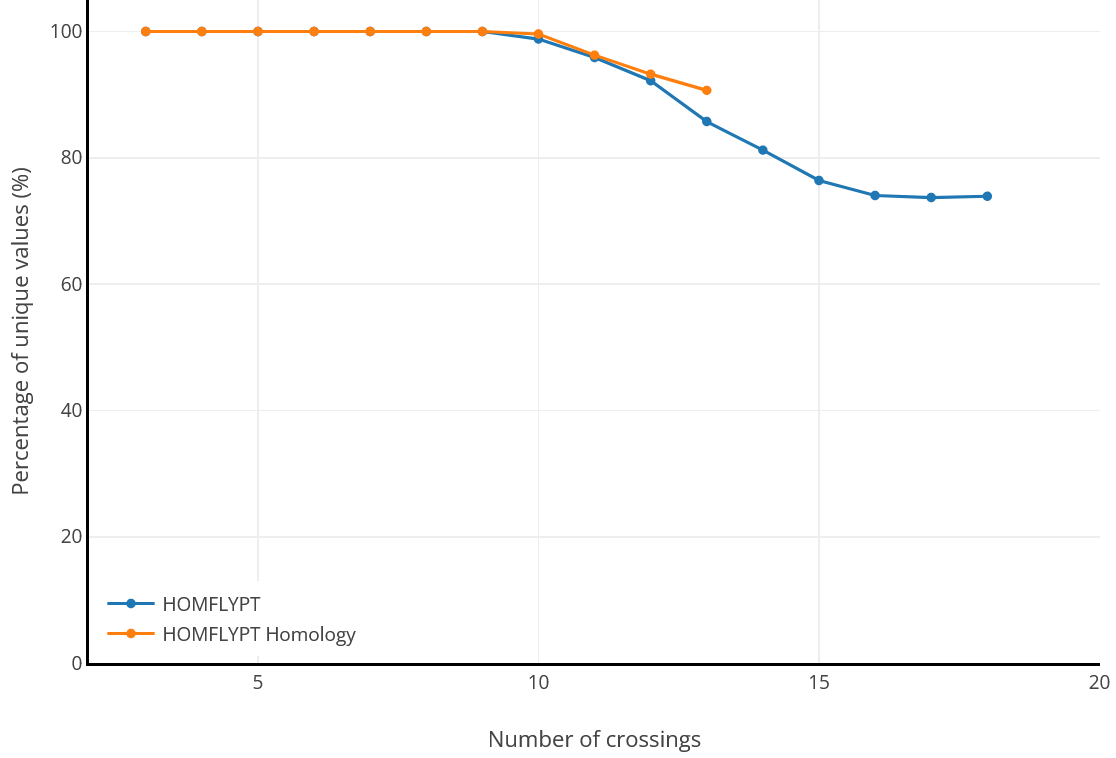}
\end{center}
\caption{Value ratios for the HOMFLYPT family. The comparison of specializations is not illustrated for the same reason as in \autoref{fig:percent-uniq-val-SL3}.}
\label{fig:percent-uniq-val-HOMFLY}
\end{figure}

The terminology on the companion website is different: its ``unique'' view gives the value ratio $v^K_Q(n)$, while its ``random'' view gives the detection ratio $d^K_Q(n)$. The four family plots \autoref{fig:percent-uniq-val-K}, \autoref{fig:percent-uniq-val-HFK}, \autoref{fig:percent-uniq-val-SL3} and \autoref{fig:percent-uniq-val-HOMFLY} display value ratios. The separate plot in \autoref{fig:jones-distribution} displays the distribution of fiber sizes, and \autoref{table:value-versus-detection} compares value and detection ratios directly.

The companion website contains the full distributions of fiber sizes, so the two statistics can be compared. For example, among the $1{,}701{,}935$ prime knots through 16 crossings, with mirrors not distinguished, the Jones polynomial takes $841{,}145$ distinct values, but only $495{,}163$ knots lie in singleton Jones fibers. Thus $v^K_J(16)=49.42\%$, while $d^K_J(16)=29.09\%$. Here the compact comparison is given in \autoref{table:value-versus-detection}.

The value ratio data are displayed in \autoref{fig:percent-uniq-val-K} to \autoref{fig:percent-uniq-val-HOMFLY} and the precise values are listed in \autoref{table:percent-uniq-val-K} to \autoref{table:percent-uniq-val-HOMFLY}. The loss of distinguishing power is already substantial and persistent in the Jones and Alexander families: through 18 crossings the value ratios have fallen to $41.61\%$ for J and $11.19\%$ for A. The stronger SL3 and HOMFLYPT polynomials level off over the presently available range, at about $74\%$, so these data do not yet give a reliable asymptotic rate for them. The companion website gives the corresponding alternating and non-alternating data, successive quotients and fiber histograms.

Value and detection ratios compress an entire fiber distribution into one number. The distribution of Jones fiber sizes through 16 crossings is shown in \autoref{fig:jones-distribution}. On log-log axes, the combined, alternating and non-alternating data are all close to straight lines over several orders of magnitude. This gives strong finite-data evidence for a heavy-tailed, approximately power-law distribution. In particular, the poor detection ratio is not caused by a handful of exceptional Jones values: nontrivial fibers occur at every visible scale, up to the largest fibers in the data.

This motivates the analog of \autoref{T:Decay}:

\begin{Conjecture}\label{C:Decay}
\textup{(}\textbf{Asymptotic failure of detection.}\textup{)}
For any fixed polynomial or homological invariant $Q$ in \autoref{Tab:Main}, there is a constant $\delta_Q\in(0,1)$ such that
\begin{gather*}
v^K_Q(n)\in O(\delta_Q^n).
\end{gather*}
Consequently, $d^K_Q(n)\in O(\delta_Q^n)$ as well.
\end{Conjecture}

This is analogous to the conclusion in \autoref{T:Decay} for prime alternating links and mutation-insensitive invariants. For general prime knots, the available range is too short to estimate $\delta_Q$ reliably for the strongest invariants, but the sustained decay in the Jones and Alexander families, the mutation data in \autoref{table:percent-jones-mut}, and the broad fiber distributions all point in the same direction. The same conjecture should hold for many further invariants and, beyond the alternating case proved here, for links.

To conclude, we highlight the observations that we find the most surprising from our data:
\begin{enumerate}

\item Categorification gives much less additional distinguishing power than one might expect. Through 16 crossings, rational Khovanov homology improves the Jones detection ratio from $29.09\%$ to $34.45\%$, a gain of only $5.36$ percentage points. At 18 crossings, its value-ratio advantage is $9.22$ percentage points. HFK is the clear exception: through 16 crossings its detection ratio is $23.66\%$, compared with $8.31\%$ for the Alexander polynomial.

\item Odd Khovanov homology performs slightly worse than standard Khovanov homology throughout the largest part of the data set: at 18 crossings their value ratios are $49.48\%$ and $50.83\%$, respectively.

\item The extra information in the strongest categorifications is tiny at the scales where it can be computed. At 13 crossings, KR3 has value ratio $86.45\%$, only $0.74$ percentage points above SL3. For HOMFLYPT homology, even the upper bound obtained by treating every missing value as unique is only $4.94$ percentage points above the HOMFLYPT polynomial at 13 crossings.

\item The SL3 and HOMFLYPT polynomials are almost indistinguishable as classifiers in the common range: at 18 crossings their value ratios are $73.84\%$ and $73.91\%$. This makes HOMFLYPT homology particularly underwhelming as a classifier relative to its computational cost; the present calculation required about 215 days and is still incomplete at 12 and 13 crossings. A direct calculation from PD presentations could improve the computational picture, but not the values already obtained.

\item The Jones fiber-size distribution is heavy-tailed. Its near-linearity on log-log axes holds separately for alternating and non-alternating knots, and shows that the failure of detection is spread throughout the data rather than concentrated in a few anomalously large fibers.

\end{enumerate}

%%%%%%%%%%%%%%%%%%%%%%%%%%%%%%%%%%%%%%%%
%%%                                  %%%
%%%            Bibliography          %%%
%%%                                  %%%
%%%%%%%%%%%%%%%%%%%%%%%%%%%%%%%%%%%%%%%%

\vspace{1em}

\begingroup
\small
\normalfont

\indent \textsc{T.K.}:\enspace {\scshape Aalto University, Department of Mathematics and Systems Analysis, Otakaari 1, 02150, Espoo, Finland,}
\href{https://orcid.org/0000-0001-9714-3988}{ORCID 0000-0001-9714-3988} \\
\indent\textit{Email address}:\enspace \texttt{tuomas.kelomaki@aalto.fi}

\vspace{0.5em}

\indent \textsc{A.L.}:\enspace {\scshape Laboratoire de Math{\'e}matiques Blaise Pascal (UMR 6620), Universit{\'e} Clermont Auvergne, Campus Universitaire des C{\'e}zeaux, 3 place Vasarely, 63178 Aubi{\`e}re Cedex, France,}\\
\href{http://www.normalesup.org/~lacabanne}{www.normalesup.org/$\sim$lacabanne}, 
\href{https://orcid.org/0000-0001-8691-3270}{ORCID 0000-0001-8691-3270}\\
\indent\textit{Email address}:\enspace \texttt{abel.lacabanne@uca.fr}

\vspace{0.5em}

\indent \textsc{D.T.}:\enspace {\scshape The University of Sydney, School of Mathematics and Statistics F07, Office Carslaw 827, NSW 2006, Australia,}
\href{http://www.dtubbenhauer.com}{www.dtubbenhauer.com}, 
\href{https://orcid.org/0000-0001-7265-5047}{ORCID 0000-0001-7265-5047}\\
\indent\textit{Email address}:\enspace \texttt{daniel.tubbenhauer@sydney.edu.au}

\vspace{0.5em}

\indent \textsc{P.V.}:\enspace {\scshape Institut de Recherche en Math{\'e}matique et Physique, 
Universit{\'e} catholique de Louvain, Che\-min du Cyclotron 2,  
1348 Louvain-la-Neuve, Belgium,}
\href{https://perso.uclouvain.be/pedro.vaz}{https://perso.uclouvain.be/pedro.vaz}, \href{https://orcid.org/0000-0001-9422-4707}{ORCID 0000-0001-9422-4707}\\
\indent\textit{Email address}:\enspace \texttt{pedro.vaz@uclouvain.be}

\vspace{0.5em}

\indent
\textsc{V.L.Z.}:\enspace {\scshape University of New South Wales (UNSW), School of Mathematics and Statistics, NSW 2052, Australia,}
\href{https://dustbringer.github.io/}{dustbringer.github.io}, \href{https://orcid.org/0009-0007-5799-6477}{ORCID 0009-0007-5799-6477}\\
\indent\textit{Email address}:\enspace \texttt{victor.l.zhang@unsw.edu.au}
\endgroup

\newpage

\begin{table}[!ht]
\fcolorbox{tomato!50}{white}{
\rowcolors{2}{spinach!25}{}
\scalebox{0.8}{\begin{tabular}{c|P{1.2cm}|P{1.2cm}|P{1.2cm}|P{1.2cm}|P{1.2cm}|P{1.2cm}}
\rowcolor{orchid!50}
n & J & K & KT1 & J+KT1 & O & K+O \\ \hline
3 & 100.0 & 100.0 & 100.0 & 100.0 & 100.0 & 100.0 \\ \hline
4 & 100.0 & 100.0 & 100.0 & 100.0 & 100.0 & 100.0 \\ \hline
5 & 100.0 & 100.0 & 100.0 & 100.0 & 100.0 & 100.0 \\ \hline
6 & 100.0 & 100.0 & 100.0 & 100.0 & 100.0 & 100.0 \\ \hline
7 & 100.0 & 100.0 & 100.0 & 100.0 & 100.0 & 100.0 \\ \hline
8 & 100.0 & 100.0 & 100.0 & 100.0 & 100.0 & 100.0 \\ \hline
9 & 100.0 & 100.0 & 100.0 & 100.0 & 100.0 & 100.0 \\ \hline
10 & 96.38 & 96.38 & 96.38 & 96.38 & 96.38 & 96.38 \\ \hline
11 & 90.13 & 91.13 & 90.76 & 91.13 & 91.13 & 91.13 \\ \hline
12 & 83.00 & 84.31 & 83.84 & 84.17 & 84.11 & 84.31 \\ \hline
13 & 73.31 & 77.50 & 77.16 & 77.44 & 74.53 & 77.50 \\ \hline
14 & 64.49 & 69.04 & 68.47 & 68.97 & 66.28 & 69.05 \\ \hline
15 & 55.74 & 60.69 & 59.85 & 60.64 & 58.40 & 60.77 \\ \hline
16 & 49.42 & 54.71 & 53.63 & 54.67 & 53.27 & 54.90 \\ \hline
17 & 44.84 & 51.93 & 50.47 & 51.90 & 50.57 & 52.33 \\ \hline
18 & 41.61 & 50.83 & 49.00 & 50.80 & 49.48 & 51.45 \\ \hline
19 & 40.17 & -     & -     & -     & -     & -     \\
\end{tabular}}
}
\caption{Value ratios for the Jones family; copyable data.}
\label{table:percent-uniq-val-K}
\end{table}

\clearpage

\begin{table}[!ht]
\centering
\begin{tabular}{c|ccc}
\rowcolor{lightergray}
$Q$ & $100v^K_Q(16)$ & $100d^K_Q(16)$ & largest fiber \\ \hline
J       & 49.42 & 29.09 & 61  \\
K       & 54.71 & 34.45 & 52  \\
KT1     & 53.63 & 33.20 & 67  \\
O       & 53.27 & 32.99 & 52  \\
A       & 18.66 &  8.31 & 923 \\
HFK2    & 38.96 & 23.66 & 410 \\
HFK2T1  & 18.15 &  8.07 & 510 \\
H       & 74.03 & 56.87 & 24
\end{tabular}
\caption{Value ratios, detection ratios and largest fiber sizes for prime knots through 16 crossings. The data on singleton fibers and fiber sizes are recorded in the online companion \cite{LaTuVaZh-quantum-big-data-code}.}
\label{table:value-versus-detection}
\end{table}

\clearpage

%\newpage

\begin{table}[!ht]
\fcolorbox{tomato!50}{white}{
\rowcolors{2}{spinach!25}{}
\scalebox{0.8}{\begin{tabular}{c|P{1.2cm}|P{1.2cm}|P{1.2cm}|P{1.8cm}}
\rowcolor{orchid!50}
n & A & HFK & HFKT1 & A+HFKT1 \\ \hline
3 & 100.0 & 100.0 & 100.0 & 100.0 \\ \hline
4 & 100.0 & 100.0 & 100.0 & 100.0 \\ \hline
5 & 100.0 & 100.0 & 100.0 & 100.0 \\ \hline
6 & 100.0 & 100.0 & 100.0 & 100.0 \\ \hline
7 & 100.0 & 100.0 & 100.0 & 100.0 \\ \hline
8 & 100.0 & 100.0 & 100.0 & 100.0 \\ \hline
9 & 94.04 & 95.23 & 92.85 & 94.04 \\ \hline
10 & 84.73 & 87.95 & 83.53 & 85.54 \\ \hline
11 & 68.78 & 79.40 & 67.16 & 70.78 \\ \hline
12 & 59.55 & 71.38 & 57.8 & 63.05 \\ \hline
13 & 43.49 & 59.31 & 42.08 & 49.61 \\ \hline
14 & 33.69 & 50.57 & 32.54 & 41.89 \\ \hline
15 & 24.55 & 43.27 & 23.82 & 35.13 \\ \hline
16 & 18.65 & 38.96 & 18.15 & 31.57 \\ \hline
17 & 14.35 & 35.82 & 14.04 & 29.61 \\ \hline
18 & 11.19 & 34.95 & 11.05 & 29.11 \\
\end{tabular}}
}
\caption{Value ratios for the Alexander family; copyable data.}
\label{table:percent-uniq-val-HFK}
\end{table}

\clearpage

\begin{table}[!ht]
\fcolorbox{tomato!50}{white}{
\rowcolors{2}{spinach!25}{}
\scalebox{0.8}{\begin{tabular}{c|P{1.2cm}|P{1.2cm}}
\rowcolor{orchid!50}
n & SL3 & KR3 \\ \hline
3 & 100.0 &  100.0 \\ \hline
4 & 100.0 &  100.0 \\ \hline
5 & 100.0 &  100.0 \\ \hline
6 & 100.0 &  100.0 \\ \hline
7 & 100.0 &  100.0 \\ \hline
8 & 100.0 &  100.0 \\ \hline
9 & 100.0 &  100.0 \\ \hline
10 & 98.79 & 98.79 \\ \hline
11 & 95.88 & 95.88 \\ \hline
12 & 92.13 & 92.54 \\ \hline
13 & 85.71 & 86.45 \\ \hline
14 & 81.20 & 82.41 \\ \hline
15 & 76.40 & 78.24 \\ \hline
16 & 74.00 & - \\ \hline
17 & 73.65 & - \\ \hline
18 & 73.84 & - \\
\end{tabular}}
}
\caption{Value ratios for the SL3 family; copyable data.}
\label{table:percent-uniq-val-SL3}
\end{table}

\clearpage

{\color{red}\textbf{Warning}:} The HOMFLYPT homology calculations for 12 and 13 crossings below are incomplete. The displayed data assumes all missing data give unique homologies, which places an upper bound on the values. There are $219/2176$ ($\approx 10.06\%$) and $3638/9988$ ($\approx 36.42\%$) missing, for 12 and 13 crossings respectively.
In both cases, knots with braid length $\geq 16$ are missing and a few knots with shorter braid length did not finish.

\begin{table}[!ht]
\fcolorbox{tomato!50}{white}{
\rowcolors{2}{spinach!25}{}
\scalebox{0.8}{\begin{tabular}{c|P{1.2cm}|P{1.22cm}}
\rowcolor{orchid!50}
n & H & HH \\ \hline
3 & 100.0 & 100.0 \\ \hline
4 & 100.0 & 100.0 \\ \hline
5 & 100.0 & 100.0 \\ \hline
6 & 100.0 & 100.0 \\ \hline
7 & 100.0 & 100.0 \\ \hline
8 & 100.0 & 100.0 \\ \hline
9 & 100.0 & 100.0 \\ \hline
10 & 98.79 & 99.59 \\ \hline
11 & 95.88 & 96.25 \\ \hline
12 & 92.20 & $\leq$93.24 \\ \hline
13 & 85.73 & $\leq$90.67 \\ \hline
14 & 81.21 & - \\ \hline
15 & 76.40 & - \\ \hline
16 & 74.02 & - \\ \hline
17 & 73.70 & - \\ \hline
18 & 73.91 & - \\
\end{tabular}}
}
\caption{Value ratios for the HOMFLYPT family; copyable data.}
\label{table:percent-uniq-val-HOMFLY}
\end{table}

\begin{table}[!ht]
\fcolorbox{tomato!50}{white}{
\rowcolors{2}{spinach!25}{}
\scalebox{0.8}{\begin{tabular}{c|P{3.2cm}|P{2.6cm}}
\rowcolor{orchid!50}
n & knots up to mutants & J \\ \hline
3  & 1      & 100.00 \\ \hline
4  & 2      & 100.00 \\ \hline
5  & 4      & 100.00 \\ \hline
6  & 7      & 100.00 \\ \hline
7  & 14     & 100.00 \\ \hline
8  & 35     & 100.00 \\ \hline
9  & 84     & 100.00 \\ \hline
10 & 249    & 93.57  \\ \hline
11 & 785    & 85.35  \\ \hline
12 & 2881   & 75.56  \\ \hline
13 & 11989  & 64.62  \\ \hline
14 & 53665  & 55.84  \\ \hline
15 & 269899 & 46.92  \\
\end{tabular}}
}
\caption{Jones value ratios on mutation classes, using Stoimenow's database \cite{Knot-database}; copyable data.}
\label{table:percent-jones-mut}
\end{table}

\end{document}